\documentclass[12pt]{article}
\usepackage{amssymb,amsmath,amsthm}
\usepackage{graphics} 
\usepackage{hyperref}
\usepackage{showlabels}

\usepackage{url}

\usepackage[numbers]{natbib}
\begin{document}

\newtheorem{theorem}{Theorem}[section]
\newtheorem{lemma}[theorem]{Lemma}
\newtheorem{condition}[theorem]{Condition}
\newtheorem{proposition}[theorem]{Proposition}
\newtheorem{remark}[theorem]{Remark}
\newtheorem{hypothesis}[theorem]{Hypothesis}
\newtheorem{corollary}[theorem]{Corollary}
\newtheorem{example}[theorem]{Example}
\newtheorem{definition}[theorem]{Definition}

\renewcommand {\theequation}{\arabic{section}.\arabic{equation}}
\def \non{{\nonumber}}
\def \hat{\widehat}
\def \tilde{\widetilde}
\def \bar{\overline}
\def\ep{\epsilon}
\def\N{{\mathbb{N}}}
\def\R{{\mathbb{R}}}
\def\Z{{\mathbb{Z}}}
\def\p{\partial}
\def\ol{\overline}
\def\l{\left}
\def\r{\right}

\date{March 29, 2017}

\title{\large{\bf Large deviations of Markov chains\\ with 
multiple time-scales}}
 \maketitle 
                                                       
\author{ \begin{center}Lea Popovic \thanks{} \\   
Department of 
Mathematics and Statistics
\\       
 lpopovic@mathstat.concordia.ca    \\
http://www.mathstat.concordia.ca/faculty/lpopovic\end{center}}

\vspace{.7in}

\begin{abstract}
\vspace{.1in}

For Markov processes evolving on multiple time-scales a combination of large component scalings and averaging of rapid fluctuations can lead to useful limits for model approximation. A general approach to proving a law of large numbers to a deterministic limit and a central limit theorem around it have already been proven in \cite{KK13} and \cite{KKP14}. We present here a general approach to proving a large deviation principle in path space for such multi-scale Markov processes.   
Motivated by models arising in systems biology, we apply these large deviation results to general chemical reaction systems which exhibit multiple time-scales, and provide explicit calculations for several relevant examples.

\vspace{.1in}


\end{abstract}

\vspace{.2in}
\tableofcontents
\vspace{.2in}



\setcounter{equation}{0}

\section{Introduction}\label{sec-introduction}
In recent years, continuous time Markov chain models 
have found extensive use in systems biology. The 
complexity of the models introduced has led to interest 
in a variety of model reduction techniques.  Some of 
these result in what are effectively laws of large 
numbers giving approximations of the model or 
subsets of the model by systems of  ordinary 
differential.  Corresponding central limit theorems for 
the deviations of the stochastic model from the 
approximating ordinary differential equation have also 
been given. In addition to the laws of large numbers and central 
limit theorems it is both natural and of some 
biological interest to consider the corresponding large 
deviation behavior of these models.

Models with what we will refer to as the ``classical'' 
scaling fit naturally into classical large deviation results 
going back to Wentzell \cite{Wen77a}, and we will review 
these briefly; 
however, our primary interest is in models with multiple time-scales.  
These models  arise from non-standard scalings of Markov chains with density-dependent rates. Dependence of the transition rates on all variables  implies a full coupling of, and an interaction between all the components. We consider arbitrary scalings of Markov chains that lead to dynamics on two dominant time-scales: a fast one - on which rapid fluctuations for a subset of components leads to geometric ergodicity; and a slow one - on which the remaining subset of components converge to a solution of a system of ordinary differential equations. 

Perhaps the simplest example in the realm of chemical reactions is a model of enzyme 
kinetics 
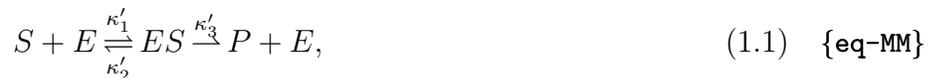
\begin{equation}\label{eq-MM}S+E\mathop{\rightleftharpoons}^{\kappa_1'}_{\kappa_2'}ES\mathop{
\rightharpoonup}^{\kappa_3'}P+E,\end{equation}
where $S$  is the substrate, $E$ the enzyme, $ES$ the 
enzyme-substrate complex, and $P$ the product.  
Under appropriate scaling of the parameters, we can 
write the model as the solution of the system
\begin{eqnarray*}
Z^N_1(t)&=&Z^N_1(0)-N^{-1}Y_1(N\int_0^t\kappa_1Z^N_1(s)Z^N_2(s)ds
)+N^{-1}Y_2(N\int_0^t\kappa_2Z^N_3(s)ds)\\
Z^N_2(t)&=&Z^N_2(0)-Y_1(N\int_0^t\kappa_1Z^N_1(s)Z^N_2(s)ds)+Y_2(
N\int_0^t\kappa_2Z^N_3(s)ds)\\
&&\qquad +Y_3(N\int_0^t\kappa_3Z^N_3(s)ds)\\
Z_3^N(t)&=&Z^N_2(0)+Y_1(N\int_0^t\kappa_1Z^N_1(s)Z^N_2(s)ds)-Y_2(
N\int_0^t\kappa_2Z^N_3(s)ds)\\
&&\qquad -Y_3(N\int_0^t\kappa_3Z^N_3(s)ds)\\
Z^N_4(t)&=&N^{-1}Y_3(N\int_0^t\kappa_3Z^N_3(s)ds),\end{eqnarray*}
where $Y_1$, $Y_2$, $Y_3$ are independent unit Poisson 
processes, and $Z^N_1, Z^N_2, Z^N_3, Z^N_4$ are the scaled amounts of substrate, free enzyme, enzyme-substrate complex, and product, respectively.  Note that $M\equiv Z_2^N(t)+Z_3^N(t)$ is constant in 
time and we will also assume independent of the scaling 
parameter $N$. The amount of substrate is an order of magnitude larger than the amount of enzymes, and hence assumed to be proportional to the scaling parameter $N$. Due to the relatively small fluctuations of the scaled amount of substrate the process $Z^N_1$ can be approximated by a deterministic one. The law of large numbers for this  system 
goes back to Darden \cite{Dard79} and is derived from the 
above 
system of equations in \cite{KK13}.  Specifically, it is 
shown that as $N\rightarrow\infty$, $Z_1^N$ converges to the solution of 
\begin{equation}\label{eq-MMlim}\dot {x}(t)=-\frac {M\kappa_1\kappa_3x(t)}{\kappa_2+\kappa_3+\kappa_
1x(t)},\end{equation}
which, of course, is simply the Michaelis-Menten 
equation. The corresponding central limit theorem for 
the scaled deviations $N^{1/3}(Z_1^N-x(t))$ is  given in 
\cite{KKP14}. 

A less straightforward example is a model of packaged virus particle production
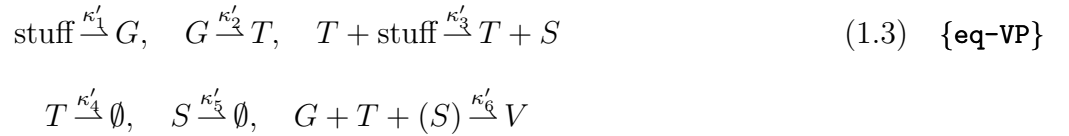
\begin{equation}\label{eq-VP}\mbox{stuff}\mathop{\rightharpoonup}^{\kappa_1'}G, \quad G\mathop{\rightharpoonup}^{\kappa_2'}T, \quad T+ \mbox{stuff}\mathop{\rightharpoonup}^{\kappa_3'}T+S\end{equation}
\[T\mathop{\rightharpoonup}^{\kappa_4'}\emptyset, \quad S\mathop{\rightharpoonup}^{\kappa_5'}\emptyset, \quad G+T+(S)\mathop{\rightharpoonup}^{\kappa_6'}V\]
where $T$ is the viral template, $G$ the viral genome, $S$ the viral structural protein that uses up resources from the cell, and $V$ is the pre-packaged material necessary for further proliferation of the virus in another cell (the structural protein is packaged, but it affects the packaging rate only in its order of magnitude). Under the appropriate scaling of the component amounts and chemical rate constants we can write the model as the solution of 
\begin{align*}
Z^N_1(t)&=Z^N_1(0)+Y_2(\int_0^t N^{2/3}\kappa_2Z^N_2(s)ds)-Y_4(\int_0^t N^{2/3}\kappa_4Z^N_1(s)ds)\\
&\qquad\qquad\qquad -Y_6(\int_0^t N^{2/3}\kappa_6Z^N_1(s)Z^N_2(s)ds)\\
Z^N_2(t)&=Z^N_2(0)+N^{-2/3}Y_1(\int_0^t N^{2/3}\kappa_1ds)-N^{-2/3}Y_2(\int_
0^t N^{2/3}\kappa_2Z^N_2(s)ds)\\
&\qquad\qquad\qquad -N^{-2/3}Y_6(\int_0^t N^{2/3}\kappa_6Z^N_1(s)Z^N_2(s)ds)\\
Z^N_3(t)&=Z^N_3(0)+N^{-1}Y_3(\int_0^t N^{5/3}\kappa_3Z^N_1(s)ds)-N^{-1}Y_5(\int_
0^t N^{5/3}\kappa_5Z^N_3(s)ds)\\
&\qquad\qquad\qquad  -N^{-1}Y_6(\int_0^t N^{2/3}\kappa_6Z^N_1(s)Z^N_2(s)ds)\\
Z^N_4(t)&=Z^N_4(0)+N^{-2/3}Y_6(\int_0^t N^{2/3}\kappa_6Z^N_1(s)Z^N_2(s)ds)
\end{align*}
where $Y_k$, $k=1,\dots,6$ are independent unit Poisson processes, and $Z^N_1, Z^N_2, Z^N_3, Z^N_4$ are the scaled amounts of template, genome, structural protein, and viral package, respectively.  
The fast fluctuating components are essentially evolving as a piecewise deterministic Markov process as defined by \cite{Dav93} with $Z_1$ a discrete component and $Z_3$ continuous.   Note that the scaled amount of template and structural protein have large fluctuations relative to their amounts, while the fluctuations of the scaled amount of genome are relatively small so the process $Z^N_2$ can be approximated by a deterministic one. 
The law of large numbers, obtained for the above system of equations by adapting the results in \cite{BKPR06}, shows that as $N\to \infty$, $Z_2^N$ converges to a solution of 
\begin{equation}\label{eq-VPlim}\dot {x}(t)=\kappa_1-\kappa_2x(t)-\kappa_6\frac{\kappa_3}{\kappa_5}\frac{\kappa_2x(t)}{\kappa_4+\kappa_6x(t)}x(t),\end{equation}
and the scaled deviations $N^{1/3}(Z_2^N-x(t))$, as can be shown by adapting the example in \cite{KKP14}, converge to a Gaussian process. In addition to two other examples, we provide large deviations for the enzyme kinetics and viral production models.\\ 
 
Our results allow a great deal of generality for the original Markov process, requiring only that it satisfy necessary technical assumptions on:\, the existence and uniqueness of limiting processes on both time-scales and some control on their exponential growth;\, geometric ergodicity of the occupation measure for the rapidly fluctuating subset of components,\,  and uniqueness for the limiting exponential operator of the remaining subset of components. They are general enough to allow the original Markov process to be a multi-scale jump-diffusion with density dependent (non-L\'evy) jump measure. Earlier results for such processes with a L\'evy measure driving the jump terms were given in \cite{KP17}. We use the same proof methodology, which relies on the general method for Markov processes  developed in \cite{FK06} based on non-linear semigroups and viscosity methods, and a generalization of Barles and Perthame limit arguments  for PDEs given in \cite{FFK12}. 

Related results for Markov chain models of chemical kinetics on two (well-separated) time-scales were recently proved in \cite{LL17} using different techniques based on the approximation and change-of-measure approach. Our results fully cover the extent of their conclusions and further extend them to more general  reaction systems with two time-scale effective dynamics.  In particular, our results allow the effective dynamics of the fast fluctuating component to be a combination of discrete and continuous variables fully linked by piecewise deterministic Markovian dynamics (PDMP) as defined in \cite{Dav93}. 
Most importantly our results do not assume the fast fluctuating component in the original Markov process to be limited to a finite or bounded state space. This is an assumption that has so far been assumed for all large deviation results on chemical reaction systems in the literature. 
The generality of the exponential weak convergence approach combined with the power of the viscosity solution technique allows all components to live in a non-compact subset of $\mathbb R^d$. 

Our large deviation principle (LDP) for general multi-scale Markov chains is given in Theorem~\ref{thm-ldp}. Our other  goal was to verify its conditions for multi-scale chemical reaction systems. Propositions~\ref{prop-binbound} and~\ref{prop-EVP} give a way to verify two technically challenging conditions: exponential compact containment, and existence of a solution to an eigenvalue problem, respectively.
We also illustrate how they apply to give the LDP in the discussed examples of enzymatic kinetics and viral production, as well as in two others.\\

{\bf Outline:} Section 2 contains the terminology for large deviations and the relevant general tools. Section 3 specifies a sequence of Conditions~\ref{cond-msconv}-\ref{cond-lyapun} that need to be verified  and the statements of the large deviation Theorem~\ref{thm-ldp} and its Corollary~\ref{cor-ldp}. Section 4 identifies specific aspects of the reaction network context that allows one to verify (or relax) the needed conditions for  multi-scale chemical reaction systems.  
Section 5 provides several examples of biologically relevant reaction systems and explicitly verifies the conditions and obtains the LDP result. The Appendix contains proofs of Theorem~\ref{thm-ldp}, Corollary~\ref{cor-ldp}, and Lemma~\ref{lem-H0cond}.



\setcounter{equation}{0}

\section{Large deviations} \label{sec-genld} For details and proofs 
regarding the following discussion see 
\cite{FK06}.  Let $(S,d)$ be a metric space.  For 
$N=1,2,\ldots$, let $X^N$ be a $S$-valued random variable.  $\{X^N\}$ 
satisfies the {\em large deviation principle\/} with {\em rate function }
$I$ if for each open set $A$,
\begin{equation}-\inf_{x\in A}I(x)\leq\liminf_{N\rightarrow\infty}\frac 
1N\log P\{X^N\in A\},\label{lb1}\end{equation}
and for each closed set $B$,
\begin{equation}\limsup_{N\rightarrow\infty}\frac 1N\log P\{X^N\in 
B\}\leq -\inf_{x\in B}I(x).\label{ub1}\end{equation}
Without loss of generality, we can require $I$ to be lower 
semicontinuous, that is $\{x:I(x)\leq c\}$ is closed for each 
$c\in {\Bbb R}$, and then we have a perhaps more 
intuitive definition of $I$,
\begin{equation}-I(x)=\lim_{\epsilon\rightarrow 0}\liminf_{N\rightarrow
\infty}\frac 1N\log P\{X^N\in B_{\epsilon}(x)\}=\lim_{\epsilon\rightarrow 
0}\limsup_{N\rightarrow\infty}\frac 1N\log P\{X^N\in\bar {B}_{\epsilon}
(x)\}.\label{icomp}\end{equation}
Typically, $\{x:I(x)\leq c\}$ will be  compact for all $c\in {\Bbb R}$. 
Then $I$ is called a {\em good rate function.}

The notion of {\em exponential tightness\/} plays the same role 
in large deviation theory that {\em tightness\/} plays in the 
theory of weak convergence. 

\begin{definition}\label{def-etight}(Exponential Tightness)
 A sequence of probability measures $\{\mu_N\}$ on $S$ 
is  {\em exponentially tight\/} if for each $a>0$, 
there exists a compact  set $K_a\subset S$ such that 
\[\limsup_{N\rightarrow\infty}\frac 1N\log\mu_N(K^c_a)\leq -a.\]
A sequence $\{X^N\}$ of $S$-valued random variables is 
exponentially tight if the corresponding sequence of 
distributions is exponentially tight.
\end{definition}  

The approach we will take to proving our large deviation 
results is based on the following theorem of Varadhan 
and Bryc (see \cite{FK06} Proposition 3.8 or \cite{DZ} Theorem 4.3.1 and 4.4.2).

\begin{theorem}\label{thm-varbryc}
Let $\{X^N\}$ be a sequence of 
$S$-valued random variables.

\begin{itemize}

\item[a)](Varadhan Lemma) Suppose that $\{X^N\}$ satisfies the 
large deviation principle with a good
rate function $I$.  Then for each $f\in C_b(S)$, 
   
\begin{equation}\lim_{N\rightarrow\infty}\frac 1N\log E[e^{Nf(X^N
)}]=\sup_{x\in S}\{f(x)-I(x)\}.\label{varadhanlemma}\end{equation}

\item[b)](Bryc formula) Suppose that the sequence $\{X^N\}$ is 
exponentially tight and  that the limit 
\begin{equation}\Lambda (f)=\lim_{N\rightarrow\infty}\frac 1N\log 
E[e^{Nf(X^N)}]\label{eq-brycformula}\end{equation}
exists for each $f\in C_b(S)$.  Then $\{X^N\}$ satisfies  the large 
deviation principle 
with good rate function 
\begin{equation}I(x)=\sup_{f\in C_b(S)}\{f(x)-\Lambda (f)\}.\label{eq-bryci}\end{equation}
\end{itemize}

\end{theorem}

We are interested in time-homogeneous Markov 
processes $\{X^N(t)\}_{t\geq 0}$ which will have sample paths in 
the Skorohod space $S=D_E[0,\infty )$. 
Assuming the limit (\ref{eq-brycformula}) exists for 
sufficiently many functions $f$, 
we can apply Theorem \ref{thm-varbryc}(b) to
\begin{equation}\Lambda_t(f,x)=\lim_{N\rightarrow\infty}\frac 1N\log 
E[e^{Nf(X^N(t))}|X^N(0)=x],\label{eq-lamlim}\end{equation}
to obtain the large deviation principle for the one 
dimensional distributions.  But if we can show
 exponential tightness for the distributions of $\{X^N\}_{t\geq 0}$ on 
$D_E[0,\infty )$, then the Markov property gives the large 
deviation principle for the finite dimensional 
distributions which in turn gives the large deviation 
principle for the processes in $D_E[0,\infty )$.

Suppose $X^N$ is a Markov processes with generator $A_N$.  Define 
\[V_N(t)f(x)=\frac 1N\log E[e^{Nf(X^N(t))}|X^N(0)=x].\]
Then by the Markov property, $\{V_N(t)\}_{t\ge 0}$ is a nonlinear 
contraction semigroup, that is,
\[V_N(t+s)f(x)=V_N(t)V_N(s)f(x),\quad s,t\geq 0\]
and
\[\sup_x|V_N(t)f_1(x)-V_N(t)f_2(x)|\leq\sup_x|f_1(x)-f_2(x)|,\]
and we can define a nonlinear {\em exponential generator\/} by
\[H_Nf(x)=\lim_{t\rightarrow 0}\frac 1t(V_N(t)f(x)-f(x))=\frac 1N
e^{-Nf(x)}A_Ne^{Nf}(x),\]
provided $e^{Nf}$ is in ${\cal D}(A_N)$.
Since it is $H_N$ that we typically 
know how to compute explicitly, it is natural to ask for 
conditions on the sequence of  generators $\{H_N\}$ that 
imply convergence of $\{V_N\}$. Observe that (\ref{eq-lamlim}) is just the convergence of 
the semigroup $V_N$.  
We define the ``limit''  as $N\to\infty$  of the sequence $H_N$ to be the set 
of $\{(f,g_{*},g^{*})\in C_b(E)\times B(E)\times B(E)\}$  for which there exists 
$f_N\in {\cal D}(H_N)$ such that for $x_N\in E$ satisfying $x_N\rightarrow x$, 
$f_N(x_N)\rightarrow f(x)$ and 
\[g_{*}(x)\leq\liminf_{n\rightarrow\infty}H_Nf_N(x_N)\leq\limsup_{
n\rightarrow\infty}H_Nf_N(x_N)\leq g^{*}(x).\] 

\

With the two examples from the introduction in mind, let us separate component-wise the notation for the multi-scale Markov process 
$Z^N=(X^N,Y^N)$ so that $X^N$ 
satisfies a law of large numbers while $Y^N$ has fast ergodic fluctuations. There are several 
complications to overcome. Since it is only $\{X^N(t)\}_{t\ge 0}$ that converges to a deterministic limit,
 we are really only interested in the 
large deviation behavior  for that sequence. Also, since the fluctuations of $\{Y^N(t)\}_{t\ge 0}$ average out, the ``limit'' of the sequence $\{H_N\}_{N\to\infty}$ will typically be a multi-valued operator.
One way to deal with identifying this limit is to select the functions $f_N(x,y)=f_0(x)+\frac 1N f_1(x,y)$ in such a way that  
$\lim_{N\to \infty}H_Nf_N(x,y)=g(x)$ for some function $g$ that is independent of $y$. For geometrically ergodic processes $Y^N$ this can typically be accomplished by solving an eigenvalue problem based on a perturbed operator for $Y^N$. However, technical challenges still remain in order to prove existence and uniqueness of the limiting semi-group by verifying the ``range condition'' for the limiting non-linear operator.

For the processes in this paper, the state space of $Z^N$ 
will always be a subset of a Euclidean space $E^N\subseteq\mathbb R^d$ 
which converges, in the sense that $E^N\subseteq E=E_X\times E_Y\subseteq \mathbb R^d$   is asymptotically dense in $\mathbb R^d$ so that for each compact $K\subset \mathbb R^d$,
\[\mathop{\lim}\limits_{N\to\infty}\mathop{\sup}\limits_{(x,y)\in E\cap K}\mathop{\inf}\limits_{(x_N,y_N)\in E^N}|(x,y)-(x_N,y_N)|=0.\]
This fact allows us to approach the problem of  convergence  of  generators $\{H_N\}$ and semi-groups $\{V_N\}$ by using the sequence of viscosity solutions of the associated Cauchy problems. Namely, for each $h\in C_b(E_X)$, the function \[u^h_N(t,x,y):=V_N(t)h(x)=\frac 1N\log E[e^{Nh(X^N(t))}|(X^N,Y^N)(0)=(x,y)]\] satisfies the non-linear partial integro-differential equation 
\begin{eqnarray}\label{eq-CauchyN}\partial_t u_N(t,x,y)&=&H_N u_N(t,x,y),\mbox{ in }(0,T]\times E_X\times E_Y; \\\ u_N(0,x)&=&h(x), \mbox{ for }(x,y)\in E_X\times E_Y.\nonumber\end{eqnarray}
The goal is to show that viscosity solutions of \eqref{eq-CauchyN} converge to a viscosity solution $u^h_0(t,x)$ of the limiting equation
\begin{eqnarray}\label{eq-Cauchylim}\partial_t u_0(t,x)&=&\bar H_0 u_0(t,x),\mbox{ in }(0,T]\times E_X; \\ u_0(0,x)&=&h(x), \mbox{ for }x\in E_X\nonumber\end{eqnarray}
where the non-linear operator $\bar H_0$ is to be identified from the limit of the non-linear generators $\{H_N\}$. In the viscosity method, existence will follow by
construction, while uniqueness will be obtained via the comparison principle. Thus, to show existence and uniqueness of the semi-group limit \eqref{eq-lamlim} using this technique, one only needs to verify the convergence of $u^h_N$ to $u^h_0$ for sufficiently many initial value functions $h$,  and to check the comparison principle for the limiting Cauchy problem \eqref{eq-Cauchylim}.

Convergence of $u^h_N$ to $u^h_0$ can be proved 
based on a general construction of subsolutions and supersolutions to two families of operators: $\{H_0(\cdot;\alpha)\}_{\alpha\in\Lambda}$ and $\{H_1(\cdot;\alpha)\}_{\alpha\in\Lambda}$, 
which are meant to ``sandwich'' the limiting operator $\bar H_0$  (see \cite{FFK12}). A comparison principle between viscosity subsolutions of $\inf_{\alpha\in\Lambda}\{H_0(\cdot;\alpha)\}$ and viscosity supersolutions of $\sup_{\alpha\in\Lambda}\{H_0(\cdot;\alpha)\}$  in conjunction with the ``operator inequality" between $\inf_{\alpha\in\Lambda}\{H_0(\cdot;\alpha)\}$ and $\sup_{\alpha\in\Lambda}\{H_0(\cdot;\alpha)\}$ will imply the desired convergence. The proof of uniqueness of the solution $u_0$ to the limiting problem then has to be shown by the weak comparison principle for the Cauchy problem  \eqref{eq-Cauchylim}. Knowledge of the eigenvalue characterization of $\bar H_0$ in the limiting problem \eqref{eq-Cauchylim} can be used in this process. The remainder of the proof of the large deviations result  comes from showing exponential tightness of $\{X^N\}$ (Definition~\ref{def-etight}) and using Bryc formula (Theorem~\ref{thm-varbryc}(b)).



\setcounter{equation}{0}
\section{Markov processes on multiple time-scales}\label{sec-mcld}

The Markov processes we are interested in have generators of the form, \begin{equation}\label{eq-A}Af(z)=\sum_k\lambda_k(z)(f(z+\zeta_k)-f(z))\end{equation}
with $k$ indexing the different jumps of size $\zeta_k$ which occur at density-dependent rates $\lambda_k(z)$. We assume the rates are non-negative, locally Lipschitz and locally bounded. 
We use powers of a parameter $N$ to scale individual component sizes $Z^N_i=N^{-\alpha_i}Z_i$ and scale density-dependent jump rates $\lambda_k(z)=N^{\beta_k}\lambda^N_k(z_N)$ leading to generators of the form, for  \mbox{$f\in\mathcal D(A_N)\subset C(E^N)$},
\begin{equation}\label{eq-AN}
A_Nf(z)=\sum_kN^{\beta_k}\lambda^N_k(z)(f(z+N^{- \underline\alpha}\zeta_k^N))-f(z))
\end{equation}
where $N^{-\underline\alpha}$ is the diagonal matrix with entries $N^{-\alpha_i}$. The nonlinear generator has the form, for \mbox{$e^{Nf}\in\mathcal D(A_N)$},
\begin{equation}\label{eq-HN}
H_Nf(z)=\frac 1N\sum_kN^{\beta_k}\lambda^N_k(z)(e^{N(f(z+N^{- \underline\alpha}\zeta_k^N)-f(z))}-1).\\
\end{equation}

\subsection{Model assumptions}
Separating the components into $Z^N=(X^N,Y^N)$ 
we need to make some assumptions on the dynamics of the rescaled Markov process on two time-scales. We start with a general look at the needed conditions.

Let $L_0$, $L_1$ be the linear operators defined on \mbox{$\mathcal D(L_0)= C_c^2(E_X)$, $\mathcal D(L_1)= C_c^2(E_X\times E_Y)$} respectively, given by 
\begin{equation}\label{eq-L0}
L_0f(z)=\sum_{k} \tilde\lambda_k(z)\tilde \zeta^X_k \cdot \nabla_x f(x)
\end{equation}
\begin{equation}\label{eq-L1}
L_1f(z)=\sum_{k:\beta_k=1} \tilde\lambda_k(z)(f(x,y+\tilde\zeta^Y_k)-f(x,y))+ \sum_{k:\beta_k>1} \tilde\lambda_k(z) \tilde\zeta^Y_k \cdot\nabla_y f(x,y)
\end{equation}  
with non-negative, locally Lipschitz and locally bounded functions $\tilde\lambda_k$ and with $\mathbb R^{|E_X|}$-valued and $\mathbb R^{|E_Y|}$-valued vectors $\tilde\zeta^X_k$ and  $\tilde\zeta^Y_k$ respectively. 
Suppose the rescaled Markov process $Z^N$ satisfies the following.

\vspace{.4cm}\noindent {\bf Condition 3.0.\, (General conditions)} \label{cond-msconv} The generator $A_N$ given by \eqref{eq-AN} satisfies
\begin{eqnarray*}\lim_{N\to\infty}\sup_{z\in E^N}&\!\! |A_Nf(z)-L_0f(z)|=0,\quad \forall f\in\mathcal D(L_0)\\
\lim_{N\to\infty}\sup_{z\in E^N}&\!\! |\frac 1N A_Nf(z)-L_1f(z)|=0,\quad\!\! \forall f\in\mathcal D(L_1)
\end{eqnarray*}
For some generator $\bar H_0$ on $\mathcal C^2_c(E_X)$ the exponential generator $H_N$ given by  \eqref{eq-HN} satisfies 
\[\lim_{N\to \infty }\sup_{z\in E^N}|H_Nf_N(z)-\bar H_0 f(z)|=0,\]
$\forall f\in\mathcal D(\bar H_0)$ and $e^{Nf_N}\in\mathcal D(A_N)$ chosen so that $\lim_{N\to \infty}\sup_{z\in E^N}|f_N(z)-f(z)|=0$.
\vspace{.4cm}

The first convergence condition essentially insures that the slow component $X^N$ has a deterministic limit. The second convergence condition insures that the fluctuations of the fast component $Y^N$ on the time-scale $tN$ have a limit that is either deterministic, a Markov chain or a piecewise-deterministic Markov chain. The latter is a Markov process with a discrete and a continuous component, where the discrete component jumps while the continuous component follows deterministic dynamics, with the jump rate and the deterministic flow determined from the value of both the discrete and the continuous component (see \cite{Dav93}). These two convergence conditions describe a separation of time scales which is a common occurrence in   stochastic models of intracellular dynamics: fast fluctuations on the time scale $tN$  are described by a limiting generator $L_1$ and slow fluctuations on the time scale $t$ are described by $L_0$  in the limit.\\ 

We next provide a sequence of explicit conditions for a large class of rescaled Markov chains with effective dynamics on two time-scales and show that the above  conditions hold. 

\begin{condition} {\bf (Scaling parameters)} \label{cond-scaling}
Let $\mathcal I_X$ and $\mathcal I_Y$ denote indices of components belonging to $X^N$ and $Y^N$ respectively. Recalling the scaling parameters $\alpha_i$ of the components $Z^N_i$ and scaling parameters $\beta_k$ of the jump rates $\lambda^N_k$, let \mbox{$\beta(i)=\max_k\{\beta_k:\zeta^N_{ik}\neq 0\}$} be the maximal jump rate that affects the component $i\in \mathcal I_X\cup\mathcal I_Y$. Then
 \[\forall i \in \mathcal I_X:\; \alpha _i\ge  1\;\mbox{ and }\; \beta(i)\le \alpha_i\]
\[\forall j \in \mathcal I_Y:\;  \alpha_j\ge 0\;\mbox{ and }\;\beta(j)\le 1+\alpha_j\]
There is at least one $i\in\mathcal I_X$ with $\alpha_i=1$ and at least one $j\in\mathcal I_Y$ with $\beta(j)=1+\alpha_j$.
\end{condition}
The first condition above requires that all the components of $X^N$ are of size at least $N$, and that jumps that change their amounts occur at rate at most $N$. The second condition insures that the separation of time-scales between fluctuations of $X^N$ and $Y^N$ is at least of order $N$. 
 The last condition implies that our choice of the scaling parameter $N$ reflects both the smallest size of the slow components and the largest separation to the time-scale of fast fluctuations. Under this condition on the scaling parameter it is reasonable to consider a large deviations result at ``speed" $N$.  

We can describe the {\it effective dynamics} on each time-scale if  we consider the {\it effective change} $(\tilde \zeta^X,\tilde \zeta^Y)$ to each component due to jumps of the process. For each $i\in\mathcal I_X$ let 
\[
\tilde\zeta^X_{ik}=\zeta^N_{ik} \mbox{ if }\beta_k=\alpha_i \mbox{ and } \tilde\zeta^X_{ik}=0\mbox{ if }\beta_k<\alpha_i,
\]
and for each $j\in\mathcal I_Y$ let 
\[
\tilde\zeta^Y_{jk}=\zeta^N_{jk} \mbox{ if }\beta_k=1+\alpha_j \mbox{ and } \tilde\zeta^Y_{jk}=0\mbox{ if }\beta_k<1+ \alpha_j.
\]
Likewise we let $\tilde \lambda_k=\lim_{N\to\infty}\lambda^N_k$ denote the {\it effective rates} of change, though in a lot of instances we will simply have $\tilde\lambda_k=\lambda_k$. Using effective changes and rates it is easy to see that Conditions~\ref{cond-scaling} imply the convergence of the sequence of generators $A_N$ from \eqref{eq-AN} to the generator $L_0$ given in \eqref{eq-L0} operating on the slow variables only, as well as the convergence of the sequence of scaled versions of these generators $\frac1N A_N$ to a generator $L_1$ given in \eqref{eq-L1}.\\

Standard assumptions on the effective processes on both time-scales are as follows. 
Let \[b_0(z)=\sum_k\tilde\lambda_k(z)\tilde\zeta^X_k, \quad b_1(z)=\sum_{k:\beta_k>1}\tilde\lambda_k(z)\tilde\zeta^Y_k, \quad c(z)=\sum_{k:\beta_k=1}\tilde\lambda_k(z)|\tilde\zeta^Y_k|.\]
\begin{condition} {\bf (Lipschitz and growth)} \label{cond-lipgro}
There exists $K_1>0$ such that $\forall z, z' \in E$ 
\[|b_0(z)-b_0(z')|+|b_1(z)-b_1(z')|+|c(z)-c(z')|\le K_1|z-z'|\]
and there exists $K_2>0$ such that $\forall z \in E$ 
\[|b_0(z)|+|b_1(z)|\le K_2|z|, \quad \sup_y c(x,y)<\infty ,\; \forall x\in E_X\]
\end{condition}
These conditions insure existence and uniqueness of the deterministic process on time-scale $t$ given by the differential operator \eqref{eq-L0}, and of the piecewise deterministic Markov process on time-scale $Nt$ given by the mixed operator \eqref{eq-L1}.  {\it However}, if existence and uniqueness of processes defined by \eqref{eq-L0} and \eqref{eq-L1} can be established by other means, for example by conditions that control their overall growth, then we can drop the above conditions and {\it assume only } that the drift coefficients $b_0,b_1$ are {\it locally Lipschitz and locally bounded} (as is implied by same assumptions on $\lambda_k$). 
We will provide conditions for dropping the global Lipschitz conditions on models of chemical reaction systems in the next Section. 

We next consider the form of the exponential generator \eqref{eq-HN} applied to functions of the form $f_N(x,y)=f(x)+\frac1Ng(x,y)$ 
\begin{eqnarray*}\label{eq-HNfN}
&&H_Nf_N(x,y)=\frac 1N\sum_kN^{\beta_k}\lambda^N_k(z)(e^{N(f_N(z+N^{- \underline\alpha}\zeta_k^N)-f_N(z))}-1)\\
&&=\sum_{k:\beta_k=1}\lambda^N_k(z)\big(e^{N(f(x+N^{-\underline \alpha}\zeta_k^N)-f(x))+(g(z+N^{- \underline\alpha}\zeta_k^N)-g(z))}-1\big)\nonumber\\
&&+\sum_{k:\beta_k>1}N^{\beta_k}\lambda^N_k(z) N^{-\underline\alpha}\zeta_k\cdot\nabla f(x)+\sum_{k:\beta_k>1}N^{\beta_k-1}\lambda^N_k(z)N^{-\underline\alpha}\zeta_k\cdot\nabla g(x,y)\nonumber\\
&&+\sum_{k:\beta_k>1}N^{\beta_k-1}\lambda^N_k(z)\big(e^{N(f(x+N^{- \underline\alpha}\zeta_k^N)-f(x))+(g(z+N^{- \underline\alpha}\zeta_k^N)-g(z))}
-N^{1-\underline\alpha}\zeta^N_k\cdot\nabla f(x)- N^{-\underline\alpha}\zeta^N_k\cdot\nabla g(x,y)-1\big)\nonumber
\end{eqnarray*}
Using 
constraints from Conditions~\ref{cond-scaling} we have that the last row above has zero limit since:\\$\beta_k+1\le 2\alpha_i$  for any $i\in\mathcal I_X$, and  $\beta_k+1=2\alpha_i$ holds only if $\alpha_i= \beta_k=1$ (as $\beta_k\le \alpha_i,\alpha_i\ge 1$); $\beta_k-1\le 2\alpha_j$ for any $j\in\mathcal I_Y$, and $\beta_k-1=2\alpha_j$  implies $\alpha_j=0,\beta_k=1$ (as $\beta_k\le \alpha_j+1$).\\ 
In the limit the effective changes $(\tilde \zeta^X,\tilde\zeta^Y)$ that we defined appear
\begin{eqnarray*}
&&\lim_{N\to\infty}H_Nf_N(x,y)
=\sum_{k:\beta_k=1} \tilde\lambda_k(z)\big(e^{\tilde\zeta^X_k \cdot \nabla f(x)}-1\big)
+\sum_{k:\beta_k=1} \tilde\lambda_k(z)e^{\tilde\zeta^X_k \cdot \nabla f(x)}\big(e^{g(x,y+\tilde\zeta^Y_k)-g(x,y)}-1\big)\\
&&\qquad\qquad\qquad\qquad +\sum_{k:\beta_k>1}\tilde\lambda_k(z) \tilde \zeta^X_k \cdot \nabla f(x)+\sum_{k:\beta_k>1}\tilde\lambda_k(z) \tilde \zeta^Y_k\cdot\nabla_yg(x,y)
\end{eqnarray*}
We let
\begin{eqnarray}\label{eq-V}
V(y;x,p)&=&\sum_{k:\beta_k=1} \tilde\lambda_k(z)\big(e^{\tilde\zeta^X_k \cdot p}-1\big)
+\sum_{k:\beta_k>1}\tilde\lambda_k(z) \tilde \zeta^X_k \cdot p
\end{eqnarray}
denote the (Hamiltonian) ``potential'' from the effective slow  process, and the operator
\begin{eqnarray}\label{eq-Lxp}
L^{x,p}_1{f(x,y)}=\sum_{k:\beta_k=1} \tilde\lambda_k(z)e^{\tilde\zeta^X_k \cdot p}{(f(x,y+\tilde\zeta^Y_k)-f(x,y))}+\sum_{k:\beta_k>1}\tilde\lambda_k(z) \tilde \zeta^Y_k\cdot \nabla_yf(x,y)
\end{eqnarray}
denote the ``perturbed'' version of the effective dynamics of the fast process. Then we 
can formulate the limit of the sequence of exponential generators on this type of functions as
\begin{eqnarray}\label{eq-H0}
\lim_{N\to\infty}H_Nf_N(x,y)=V(y;x,\nabla_xf(x))+e^{-g(x,y)}L_1^{x,\nabla_xf(x)}e^{g(x,y)}.
\end{eqnarray}
Since the limiting operator appears multi-valued, we aim to select $g$ in such a way that the right-hand side does not depend on $y$. 
In other words, we look for $g$ so that the limiting operator in \eqref{eq-H0} using $p=\nabla_xf(x)$ 
solves the eigenvalue problem:  for all $x\in E_X, p\in \mathbb R$
\begin{equation}\label{eq-EVP}
\left(V(y;x,p)+L_1^{x,p}\right)e^{g(x,y)}=\bar H_0(x,p)e^{g(x,y)}.
\end{equation}
In that case $\bar H_0(x,p)$ does not depend on $g$ and we have \[\lim_{N\to\infty}H_Nf_N(x,y):=\bar H_0(x,\nabla_xf(x)).\]
If there exists an eigenfunction $g$ such that the eigenvalue problem \eqref{eq-EVP} holds, then  Conditions~\ref{cond-scaling} provide a definition via the equation \eqref{eq-EVP} of a generator $\bar H_0$ operating on the slow variables only. We will use the PIDE method to show the sequence of exponential generators $H_N$ when applied to functions of the form $f_N$ converge to $\bar H_0$ applied to $f$. \\

Having established a candidate for the limiting operator $\bar H_0$ one further needs to establish exponential tightness for the slow process $\{X^N\}$ and for the occupation measures of the fast process  $\{\Gamma^N(\cdot,C)\}$ where 
\[\Gamma^N(t,C)=\int_0^t \mathbf 1_{C}(Y^N(s))ds.\]
The next two conditions insure multiplicative ergodicity of the occupation measures and exponential stability of the dynamics for the effective fast process perturbed by the directional changes of the slow process. We will show that exponential tightness will follow when we combine these conditions with control of the growth of the slow process and with convergence of the exponential generators. 

\begin{condition} {\bf (Transition density)} \label{cond-ergod}
For each $x\in E_X, p\in \mathbb R$ the process $Y^{x,p}$ defined by the generator \eqref{eq-Lxp} is Feller continuous with transition probability $p^{x,p}_t(y,dy)$ which at $t=1$ has a positive density with respect to some reference measure $\alpha(dy)$ on $E_Y$. 
\end{condition}

\begin{condition} {\bf (Lyapunov I)} \label{cond-lyapun}
There exists a positive function $\varphi(\cdot)\in C^1(E_Y)$ with compact level sets and such that for each compact $K\subset \mathbb R$, $\theta\in(0,1]$ and $l\in \mathbb R$ there exists a compact set $A_{l,\theta,K}\subset E_Y$ satisfying
\begin{equation}\label{eq-Lyapunov}
\{y\in E_Y:   -\theta e^{-\varphi}L_1^{x,p}e^{\varphi}(y)-|V(y;x,p)|\leq l\}\subset A_{l,\theta,K}, \quad \forall p\in K, \forall x\in E_X,
\end{equation}
and for each $x\in E_X, p\in\mathbb R^{|E_X|}$ there exists $K_{x,p}>-\infty$ such that $\forall z \in E$ 
\[V(y;x,p)\ge K_{x,p},\; \forall y\in E_Y\]
\end{condition}

The transition density condition was given by \cite{DV83} to insure a large deviation result for the occupation measures of a general Markov process $Y^{x,p}$. Together with the second condition it immediately implies that for each $(x,p)$ there exists a stationary distribution for $Y^{x,p}$ (see \cite{FK06} Lemma ~11.23). Different versions of such conditions are given in \cite{FK06} Condition 11.21  (also, see Appendix B for other related references). We will discuss conditions for piecewise deterministic Markov processes given by a perturbation of \eqref{eq-L1} which can be used to verify multiplicative ergodicity and exponential stability in the next Section.
 

In some cases it will suffice to verify exponential stability of the effective fast process using an $(x,p)$-dependent function for the Lyapunov condition. This condition together with a version of a condition on the transition kernels implies Condition (DV3+) in \cite{KM05} for verifying multiplicative ergodicity and establishing a large deviation principle for the the occupation measure of single time scale Markov processes.

\begin{condition} {\bf (Lyapunov  II)} \label{cond-multiplyapun}
For each $p$ in compact set $K\subset \mathbb R$ and for each $x\in E_X$ there is a positive function $\varphi_{x,p}(\cdot)\in C^1(E_Y)$ with compact level sets such that there exist $c>1$ and $d_{x,p}<\infty$ satisfying\\
\begin{equation}\label{eq-multiplyapun}
e^{-\varphi_{x,p}}L_1^{x,p}e^{\varphi_{x,p}}(y)\le -c |V(y;x,p)|+d_{x,p}, \quad \forall y\in E_Y. 
\end{equation}
where we assume that $|V(y;x,p)|$ has compact  level sets.
\end{condition}

Finally, uniqueness of the limiting operator $\bar H_0$ needs to be either assumed or established. We will consider the first option in our main (upcoming) Theorem and the second in its Corollary.

\subsection{Large deviation principle for the two time-scale model}

We can now present the large deviation result for the Markov process on two time-scales.
\begin{theorem}\label{thm-ldp} 
Assume Conditions~\ref{cond-scaling}-\ref{cond-lyapun} hold as well as the weak comparison principle for the Cauchy problem:
\begin{eqnarray}\label{eq-Cauchy0}
\partial_t u_0(t,x)&=& \bar H_0(x,\partial_x u_0(t,x)),\; \mbox{ for } (t,x)\in(0,T]\times E_X\\
u_0(0,x)&=&f(x),\;\mbox{ for } x\in E_X\nonumber
\end{eqnarray}
 with $\bar H_0$ defined by \eqref{eq-H0}. Then the sequence $\{X^N(t)\}$ is exponentially tight and satisfies a large deviation principle with speed $N$ and good rate function $I$ given by the variational principle:
\begin{equation}\label{eq-rate}
I(x,x_0,t)=\sup_{f\in C_b(E_X)}\{f(x)-u_0^f(t,x_0)\}
\end{equation}
where $u_0^f$ is the unique continuous viscosity solution of \eqref{eq-Cauchy0}.
\end{theorem}

The proof of Theorem~\ref{thm-ldp} is predominately following the arguments in \cite{KP17}. We show  that this argument can be carried through for processes generated by more general (non-L\'evy) jump terms and without diffusions. Two key ingredients make the proof work:  multiplicative ergodicity of the effective fast process (perturbed by the value of the slow process and its drift) which insures exponential control on the behaviour of its occupation measure; and exponential stability for the dynamics of the effective fast process (in a potential generated by the slow process) which insures exponential compact  containment and tightness. The proof is given in Section~\ref{subsec-proofthm} in the Appendix.

Much of the proof concentrates on proving the convergence of solutions to Cauchy problem \eqref{eq-CauchyN} via a sequence of sub- and super-operators which sandwich the limiting operator $\bar H_0$. A part of this proof can be significantly simplified when we know the principle eigenvalue problem with $\bar H_0$ has a solution and we make use of its associated eigenfunction. Being able to solve the eigenvalue problem explicitly relies on the form of the rates $\lambda_k$ as well as on some of the structural properties of the Markov process itself. For polynomial $\lambda_k$ there is a class of models (see Condition~\ref{cond-conversion}) for which we can explicitly solve for the eigenvalue and eigenfunction of \eqref{eq-EVP}. We will discuss a procedure for solving $\bar H_0$  for models of chemical reaction systems in the next Section.

Having an expression for $\bar H_0(x,p)$ also allows for a relatively easy verification of the comparison principle for the limiting Cauchy problem \eqref{eq-Cauchy0} by using the following result. The proof of this result is based on the theory of discontinuous viscosity solutions (see either \cite{BD} Chapter V, or \cite{FS06} Chapter VII) and is contained in Section~\ref{subsec-prooflem} of the Appendix.  For definitions of lower semicontinuous viscosity sub- and upper semicontinuous super-solutions for discontinuous functions and of the weak comparison principle we refer to \cite{FS06} VII Def 4.2 and VII Def 7.1. 

\begin{lemma}\label{lem-H0cond} 
Suppose $u_1$ and $u_2$ are,respectively, a bounded upper semicontinuous (USC) viscosity sub-solution and a bounded lower semicontinuous (LSC) viscosity  super-solution of  \eqref{eq-Cauchy0}
for some $T>0$ and  $E_X\subset \R^d$. Either of the following conditions are sufficient for the weak comparison principle for \eqref{eq-Cauchy0} to hold:
\begin{itemize}
\item[(a)]  $\ol H_0$ is such that for all $\lambda\ge 1, R>0$ and for all $|p|,|q|\le 1$, $|x|,|y|\le R$ and for some continuous non-decreasing functions $\omega_R, \,\tilde \omega_1:\R_+\mapsto\R_+$ with $\omega_R(0)=\tilde \omega_1(0)=0$
\begin{equation}\label{eq:H2}\begin{aligned}
\hspace{-5mm}\ol H_0&(y, \lambda(x-y)+p)-\ol H_0(x,\lambda(x-y)+q)\le \omega_R(|x-y|+\lambda |x-y|^2)+\tilde\omega_1(|p-q|)\end{aligned}
\end{equation}
\item[(b)]  $\ol H_0$ is such that for all $\lambda\ge 1, R,\ell>0$ and for all $|p|,|q|\le 1$, $|x|,|y|\le R$ \underline{with $\lambda|x-y|<\ell$}, and for continuous non-decreasing functions $\omega_{R,\ell}, \,\tilde \omega_1:\R_+\mapsto\R_+$ with \mbox{$\omega_{R,\ell}(0)=\tilde \omega_1(0)=0$}  the inequality \eqref{eq:H2} holds;\, 
and 
$\ol H_0$ satisfies a coercivity condition in $p$
\begin{equation}\label{eq:coercive}
\frac{\ol H_0(x,p)}{|p|}\to c_x>1\mbox{ as }|p|\to\infty \mbox{ uniformly with respect to $x$},
\end{equation} 
\end{itemize}
\end{lemma}

We now have a more easily verified large deviation result for multi-scale Markov processes.

\begin{corollary}\label{cor-ldp}
Assume Conditions~\ref{cond-scaling}-\ref{cond-ergod} and Condition~\ref{cond-multiplyapun} hold, without necessarily a lower bound on $V(y;x,p)$ but assuming $|V(y;x,p)|$ has compact level sets. Suppose the eigenvalue problem \eqref{eq-EVP} for  $\bar H_0(x,p)$ can be explicitly solved, and either of the conditions (a) or (b) from Lemma~\ref{lem-H0cond} hold. Then the large deviation principle  for $\{X^N\}$ as stated in Theorem~\ref{thm-ldp} holds.
\end{corollary}

The proof of Corollary~\ref{cor-ldp} is centered on simplifying the construction of the approximating sequence of operators from the proof of Theorem~\ref{thm-ldp}. In order to do so it uses the eigenfunction in the principle eigenvalue problem for $\bar H_0$. The proof is given in Section~\ref{subsec-proofcor} in the Appendix. Its practical advantages are that it allows one to use a Lyapunov function dependent on the variable for the slow process; and it simplifies establishing uniqueness of the limiting operator $\bar H_0$. 
Note that if the state space for the fast process is compact one can trivially take $\varphi=0$ in either \eqref{eq-Lyapunov} or \eqref{eq-multiplyapun}. Moreover, if the state space is non-compact, the process of finding the explicit solution for the eigenvalue problem can be adapted to find the Lyapunov function $\varphi_{x,p}$ as well. We will discuss finding candidates for the Lyapunov functions in the next Section.

To extend the above results, from either Theorem~\ref{thm-ldp} or Corollary~\ref{cor-ldp}, from finite dimensional distributions of $\{X_N(t)\}$ to pathwise large deviation principle for $\{X^N\}$ can be done by the argument from Theorem 4.28 of \cite{FK06}. This will require using a variational representation of the operator $\bar H_0(x,p)$ and characterizing an expression for the rate function from its Fenchel-Legendre transform (\cite{FK06} Section 8.6.1, \cite{F99}). This is the characterization of the rate functions given in \cite{LL17}.\\


\def\ep{\epsilon}
\def\ol{\overline}
\def\ul{\underline}
\def\l{\left}
\def\r{\right}
\def\p{\partial}

\setcounter{equation}{0}
\section{Chemical reaction networks on multiple time-scales}\label{sec-crld}

We now apply our large deviation results to models of multi-scale chemical reaction networks. 
In systems biology continuous-time Markov chains have found an important application for modelling chemical reactions describing cellular metabolic, gene regulatory and signal transduction processes.  
A stochastic model of an intra-cellular chemical reaction network treats the system as a continuous time Markov chain with generator of the form \eqref{eq-A} whose state $Z$ is a vector giving the number of molecules of different types of chemical species that are relevant. Each reaction is modeled as a possible transition for the state. The model for the $k$th reaction, for each $k$, is determined by a vector of inputs $\nu'_k$ specifying the number of molecules of each chemical species that are consumed in the reaction, and a vector of outputs $\nu_k$ specifying the number of molecules of each species that are produced in the reaction. Transition rate for the $k$th reaction $\lambda_k(z)$ is a function of the state $z$, and the state change of $k$th reaction is given by $\zeta_k=\nu'_k-\nu_k$. 

Reaction rates in chemical networks are most commonly modelled by dynamics of {\it mass-action} type: in the stochastic version of the law of mass action, the rate function is proportional to the number of ways of  selecting the molecules that are consumed in the reaction:
\begin{equation*}\label{eq-massact}\lambda_k(z)=\kappa'_k\prod_i\nu_{ik}!\prod_i{{z_i}\choose {\nu_{
ik}}}=\kappa'_k\prod_iz_i(z_i-1)\cdots (z_i-\nu_{ik}+1).\end{equation*}
making $\lambda_k(z)$ a product of abundances of all the species going into the reaction and of a {\it chemical reaction constant} $\kappa'$.
Physically, $|\nu_k|=\sum_i\nu_{ik}$ is usually assumed to be less than or equal to two, which makes the task of controlling the growth of coordinates much easier.
In some models a few reaction rates can be given in terms of a sigmoid function. This is often the result of approximating the overall outcome of a subnetwork of reactions through a model reduction procedure. Nonetheless, in all cases encountered in the literature the jump rates are locally Lipschitz and locally bounded, though, as is the case in binary reactions, they are not always globally Lipschitz. \\

We next discuss satisfying the conditions needed for Theorem~\ref{thm-ldp} and for Corollary~\ref{cor-ldp} to apply in models of chemical reaction networks, and specify when they can be relaxed.
Multiple scalings of intra-cellular chemical reaction processes arise naturally due to low copy numbers of various key chemical species types. In other words, species $i$ can be in abundance of order $N^{\alpha_i}$ in the system, where typically $N$ is the order of magnitude of the most abundant species and $\alpha_i\in[0,1]$. Each species is then represented by a component with its rescaled size $Z_i^N=N^{-\alpha_i}Z_i$. In addition, varying strengths in chemical bonds formed or broken by different  reactions lead to different orders of magnitude, that can be expressed  in terms of $N$, for the different chemical reaction constants $\kappa_i'$.   We combine the scaling of  constants $\kappa_i'$ together with the effect of rescaling species amounts into a single rescaling of the reaction rates as $\lambda_k(z)=N^{\beta_k}\lambda_k(z^N)$, where typically $\beta_k>0$. The rescaled chemical reaction network then becomes a multi-scale Markov process on $E^N\subset \mathbb R^d_+$ with generator  of the form \eqref{eq-AN}.

We are interested in models of intercellular reactions whose dynamics has a multi-scale behaviour with two separated time-scales as in {\bf Conditions~\ref{cond-scaling}}. This means that in the reaction networks there is a group of \emph{discrete} species is present in small counts $O(1)$ and the rest, \emph{continuous} species, are present in larger amounts approximated by $O(N)$. Their intertwined reactions dynamics consists of fast $O(N)$ fluctuations for the discrete species and slow $O(1)$ changes of the continuous species. The effective changes $\tilde \zeta_k$ are typically the same as the net reaction changes $\zeta_k$, except in instances when the order ${\alpha_i}$ of some species $i$ with $\zeta_{ki} \neq 0$ is smaller than the order ${\beta_k}$ of $\lambda_k$ and the effective change $\tilde \zeta _{ki}=0$ (see, for example, species $V$ in reaction (6) of the viral production model from the Introduction). 

A reaction network often involves a few binary reactions with a quadratic rate, so for the effective processes to exist and be unique one needs  the overall depletion rates for each species to balance out their overall production rates, insuring that all amounts are globally stable. Many multi-scale Markov models have global existence and uniqueness with only local Lipschitz and growth of jump rates: the deterministic dynamics of the effective slow process (governed by growth of the drift $b_0(z)$) may be globally controlled; and the piecewise deterministic Markov process (which is often simply a Markov chain) may have global stability (determined by the local drift $b_1(z)$ of the continuous part and the overall jump rates $c(z)$ of the discrete part) for each fixed value of the slow process. In that case {\bf Condition~\ref{cond-lipgro}} is not necessary and instead {\bf Condition~\ref{cond-binbound}} on the binary rates (discussed in the next subsection) should be verified. We prove in Proposition~\ref{prop-binbound} a truncation argument using which the requirement of globally Lipschitz coefficients can be relaxed (both the enzyme kinetics and the viral production model have some binary reaction rates that can be accommodated by this result).

The effective fast process perturbed by the direction of change of the slow process (given by generator \eqref{eq-Lxp}) is in general a piecewise deterministic Markov process (though sometimes it is simply a Markov chain). Exponential stability and positivity of its transition density assumed in {\bf Condition~\ref{cond-ergod}} can be verified using recent results on piecewise deterministic processes \cite{BH12} (see also \cite{CH15} Ch 4). To summarize the conclusions, suppose the discrete component of this PDMP  has jump rates such that their infimum over the continuous component yields an irreducible and positive recurrent Markov chain. Suppose that for all  values of the discrete component the generator of the continuous component satisfies a drift condition with respect to the same Lyapunov function. Furthermore suppose that the flow satisfies a Hoermander-type bracket (hypoellipticity)  condition at a point of the continuous component. Then,  at this point, the process has a positive transition probability with non-trivial absolutely continuous part with respect to  Lebesgue measure. 
These results are particularly useful when the state space for the effective fast process is non-compact and having a positive transition density is non-trivial.

In models where the state space $E_Y$ for the fast process is compact verifying the exponential Lyapunov {\bf Condition~\ref{cond-lyapun}} is unnecessary. In the non-compact case we need to look for a candidate function $\varphi$ based on both the dynamics of perturbed effective fast process with generator $L_1^{x,p}$, and on the potential function $|V(y;x,p)|$. For multi-scale processes for which we can explicitly solve the eigenvalue problem \eqref{eq-EVP} we can relax this condition and instead verify the  Lyapunov {\bf Condition~\ref{cond-multiplyapun}}, and look for a candidate function $\varphi_{x,p}$ dependent on the value and direction of change of the slow process. 
Moreover, such a function can be found by a similar procedure according to which one can find the eigenfunction of the EVP (see, for example, the models of down-regulation and of viral production in Examples section, both of which have an unbounded fast variable).  For practical purposes, with {\bf Condition~\ref{cond-conversion}} and Proposition~\ref{prop-EVP} we provide a procedure for solving the EVP explicitly.

 \subsection{Truncating the jump rates}\label{subsec-binbound}

Many chemical reaction networks involve a few binary reactions whose rates have quadratic growth. Our goal is to truncate these jump rates and prove the large deviation result using Theorem~\ref{thm-ldp} on a process with  truncated rates. Let us illustrate this on the examples from the Introduction. 
In the model of enzymatic kinetics \eqref{eq-MM} the sum of enzyme and enzyme-substrate abundances together is conserved by the system. Since $Z_2(t)+Z_3(t)\equiv M$ is constant in time, this means that the rate of enzyme-substrate production $S+E\rightharpoonup ES$ satisfies $\lambda_1(z)=\kappa_1z_1z_2\le \kappa_1 M z_1$, and so all the rates in the system are at most linear. 
In the model of viral production \eqref{eq-VP}, there is a quadratic reaction rate $\lambda_6(z)=\kappa_6z_1z_2$ for the viral packaging $G+T+(S)\rightharpoonup V$.  
There is no conservation law in the system, but the viral genome G is slowly varying and its rate of increase is given only by $\emptyset\rightharpoonup G$ and is bounded by a constant. This allows us to have exponential control on how large the factor $G$ can get. As we will see, this  will imply that from the point of view of large deviations all reaction rates can be replaced with versions that are at most linear.  

We make the following assumption on the jump rates.
\begin{condition}\label{cond-binbound}
Suppose that for each $k$ there are non-negative constants $\theta_{k,0}$ and $\theta_k=\{\theta_{k,i}\}$ so that
\begin{equation}\label{eq-binbound}
\lambda_k(z)\leq (\theta_{k,0}+\sum_{i}\theta_{k,i}z_i)(\sum_iz_i).\end{equation}
and for each $k$ we have either:

 {\bf (i)} all species $i$ for which $\theta_{k,i}>0$ are part of a general conservation law, that is, $\theta_{k,i}>0$  implies $\tilde\theta_i>0$, where ${\tilde\theta}=\{\tilde\theta_i\}$ is such that all net changes satisfy
\[\zeta_k\cdot{\tilde\theta}=0, \;\; \forall k;\vspace{-3mm}
\]
or: 

{\bf (ii)} all species $i$ for which $\theta_{k,i}>0$ are part of the slow process, and rates of all those reactions with a net increase of species combinations given by $\theta_k$ satisfy
\[\sum_{k:\zeta_k\cdot \theta_{k,i}>0} \lambda_k(z)\le C', \; \mbox{ for some }C'<\infty.\]
\end{condition}
For rates that are specified in ``mass-action" form the Condition~\eqref{eq-binbound} requires that in all  binary reactions at least one of the reactants $S_i$ must satisfy $\theta_{k,i}>0$. Also, if that reactant is part of the fast process it must be part of a conservation law (there can be multiple conservation laws in the system), or if it is part of the slow process it is created only by reactions with bounded rates, for example, reactions such as $\emptyset\rightharpoonup S_i$ or reactions $S_j\rightharpoonup S_i$ for some species $S_j$ of bounded abundance.

We need to introduce a notion that will describe what happens to the process once we truncate its jump rates, and justify why it is useful in the context of large deviation results (this is Theorem 4.2.16 in  \cite{DZ}).

\begin{definition} {\it (Exponential Approximation)}
The sequence $\{X^{N,M}\}$ is an {\rm exponentially good approximation} of  $\{X^N\}$ if for every $\epsilon>0$ 
\[\lim_{M\to\infty}\limsup_{N\to\infty}\frac 1N \log P\{|X^N-X^{N,M}|>\epsilon\}=-\infty.\]
 \end{definition}
 
 \begin{theorem}
 Suppose for each $M$ the sequence $\{X^{N,M}\}$ satisfies a large deviation principle with rate function $I_M$ and suppose $\{X^{N,M}\}$ is an exponentially good approximation of $\{X^N\}$. Then
 
 a) $\{X^N\}$ satisfies a {\rm weak} large deviation principle (meaning that \eqref{lb1} holds for each open set $A$ while \eqref{ub1} only holds for each compact set $B$) with rate function 
 \[I(y)=\sup_{\epsilon>0}\liminf_{M\to\infty}\inf_{z\in B_{\epsilon}(y)} I_M(x).\]

b) If $I(\cdot)$ is a {\rm good} rate function (sublevel sets are compact) 
and for each closed set $B$
\[\inf_{y\in B}I(y)\le \limsup_{M\to\infty}\inf_{y\in B}I_M(y),\]
then $\{X^N\}$ satisfies the large deviation principle with rate function $I$.
 \end{theorem}
 
 The following result is very useful in some models of chemical reaction systems.
\begin{proposition}\label{prop-binbound} Assume Condition~\ref{cond-binbound} holds and that the initial value satisfies $P\{X^N(0)\cdot \theta_k\le C\}=1$ for some $C<\infty$ and all $\theta_k$ as in \eqref{eq-binbound}. Then for all large enough $M<\infty$ replacing the jump rates $\lambda_k(z)$ of the Markov process $X^N$ by $\lambda_k(z)\wedge(\theta_{k,0}+M)(\sum_iz_i)$ produces a sequence of processes $\{X^{N,M}\}$ that is an exponentially good approximation of $\{X^N\}$, that is
\[\limsup_{N\to\infty}\frac1N \log P\{\sup_{s\le t}|X^N(t)-X^{N,M}(t)|>0\}\le c(t)-M.\]
\end{proposition}
\begin{proof}
In case (i) holds for reaction $k$ the assumption on $\tilde\theta$ implies that $\tilde\theta_{i}=a_{k,i}\theta_{k,i}$ for some $a_{k,i}>0$. The conservation law implies that $X^N(t)\cdot \tilde \theta=X^N(0)\cdot \tilde \theta$  for $t>0$. The assumption on the initial value implies 
\[X^N(t)\cdot \theta_k \le (\max_{k,i} \frac{1}{a_{k,i}}) X^N(t)\cdot \tilde \theta=(\max_{k,i} \frac{1}{a_{k,i}}) X^N(0)\cdot \tilde \theta\le \frac{\max_{k,i} {a_{k,i}}}{\min_{k,i}{a_{k,i}}}\, C.\]
So as soon as $M\ge \frac{\max_{k,i} {a_{k,i}}}{\min_{k,i}{a_{k,i}}}\, C$ we have that $P\{\sum_i \theta_{k,i}X^N_i(t)\le M\}=1$, and we can replace the rates $\lambda_k(z)$ by $\lambda_k(z)\wedge (\theta_{k,0}+M)(\sum_iz_i)$ without altering the process.

In case (ii) holds we make use of the martingales
\[\exp\Big\{N\big(X^N(t)\cdot\theta_k-X^N(0)\cdot\theta_k-\int_0^t\sum_k N^{\beta_k-1}\lambda^N_k(X^N(s))\big(e^{N^{1-\underline\alpha}\zeta_k\cdot\theta_k}-1)\big)ds\big)\Big\}\]
to obtain for any $M<\infty$ and stopping time $\tau^N$ for $X^N$ that
\begin{eqnarray*}
&&\limsup_{N\rightarrow\infty}\frac 1N\log P\{\sup_{s\leq t}X^N(s)\cdot\theta_k\geq M\}\\
&&\leq\limsup_{N\rightarrow\infty}\frac 1N\log\frac {E[e^{N X^N(t\wedge\tau^N)\cdot\theta_k}]}{e^{NM}}\\
&&\leq\limsup_{N\rightarrow\infty}\frac 1N\log\frac {E[\exp\{N \big(X^N(0)\cdot\theta_k+\int_0^{t\wedge\tau^N}\lambda_k(X^N(s))\sum_{k:\zeta_k\cdot\theta_k>0} N^{\beta_k-1}(e^{N^{1-\underline\alpha}\zeta_k\cdot\theta_k}-1\big)ds)\}]}{e^{NM}}\\
&&\leq C+C't(e^{\max_{k:\beta_k=1,\; \zeta_k\cdot\theta_k>0}\tilde\zeta^X_k\cdot\theta_k}+\max_{k:\beta_k>1,\zeta_k\cdot\theta_k>0}\tilde\zeta^X_k\cdot\theta_k)-M=c(t)-M\end{eqnarray*}
in the last inequality we used the assumption $P\{\sum_{k:\zeta_k\cdot\theta_k>0}\lambda_k(X^N(t))\leq C'\}=1$ for $\forall t\ge 0$, and that since $\zeta_k\cdot\theta_k=\zeta^X_k\cdot\theta_k$
\begin{eqnarray*}&&\lim_{N\rightarrow\infty}\sum_{k:\zeta_k\cdot\theta_k>0} N^{\beta_k-1}(e^{N^{1-\underline\alpha}\zeta_k\cdot\theta}-1)=\sum_{k:\beta_k=1,\; \zeta_k\cdot\theta>0} (e^{\tilde\zeta^X_k\cdot\theta}-1)+\sum_{k:\beta_k>1,\;\zeta_k\cdot\theta>0}\tilde\zeta^X_k\cdot\theta\\
\end{eqnarray*}
Let $X^{N,M}$ be the process with the same chemical reaction network properties (same $\zeta_k, \alpha_i, \beta_k, Y_k$) except we replace the rates $\lambda_k(z)$ by $\lambda_k(z)\wedge (\theta_{k,0}+M)(\sum_iz_i)$. 
On the event ${\{\sup_{s\leq t}X^N(s)\cdot\theta_k< M\}}$ the processes $X^N$ and $X^{N,M}$ are the same, hence
\[\limsup_{N\rightarrow\infty}\frac 1N\log P\{\sup_{s\leq t}|X^N(s)-X^{N,M}(s)|>0\}\leq\limsup_{N\rightarrow\infty}\frac 1N\log P\{\sup_{s\leq t}X^N(s)\cdot\theta_k\geq M\}\le c(t)-M\].
\end{proof}

\subsection{Solving the eigenvalue problem}\label{subsec-eigenvprob}

In order to apply Theorem~\ref{thm-ldp} we also need to show the comparison principle for the limiting Cauchy problem \eqref{eq-Cauchy0}. It insures that we can identify the limiting semi-group \eqref{eq-lamlim} in Bryc formula from our convergence arguments. Theory for uniqueness of viscosity solutions says that the comparison principle will hold if the conditions provided by Lemma~\ref{lem-H0cond} are satisfied. 
 They are similar to the condition for proving comparison principle for operator themselves (see Lemma 9.2 \cite{FK06} or Proposition 7.7 \cite{CIL92}). But the advantage of the convergence approach via solutions of the Cauchy problem (where one needs to use only bounded solutions), rather than the operator, is that we only need to check these conditions for variables in a compact subset of values for $x$ and $p$ if we have an operator that is coercive in $p$ uniformly over $x$. This is most easily verified if $\bar H_0$ can be explicitly calculated. 

In order to explicitly calculate $\ol H_0$ we need to solve the eigenvalue problem \eqref{eq-EVP}. 
We first make the following assumptions on the jump rates in case the effective dynamics of the fast process is a pure Markov chain.

\begin{condition}\label{cond-conversion}
Suppose that $\tilde\zeta^Y_{k,i} \in\{-1,0,1\}\; \forall k,\forall i$; there are constants $\theta_{k,0},\{\theta_{k,i}\}$ $\in \{0,1\}$ satisfying $\theta_{k,0}+\sum_i\theta_{k,i}=1$ such that
 \[\lambda_k(x,y) = \lambda_k(x)(\theta_{k,0}+\sum_i\theta_{k,i} y_i);\]
and for each species $j$ for which $\exists k,i$ such that $\theta_{k,i}= 1$ and $\tilde\zeta^Y_{k,j}=-1$   we have that either:

{\bf (i)} a species $j'\neq j$ for which $\exists k'$ (possibly $k$) such that $\theta_{k',i}=1$  and $\tilde\zeta^Y_{k',j'}\neq 0$ is part of a conservation law together with species $j$, that is, both $\tilde \theta_j>0 \mbox{ and } \tilde \theta_{j'}>0$ where  $\tilde \theta=\{\tilde \theta_i\} $ is such that effective net changes for the fast process satisfy
\[\tilde\zeta_k^Y\cdot \tilde\theta =0, \;\; \forall k;\vspace{-3mm}\]
or: 

{\bf (ii)} a species $j'\neq j$  for which $\exists k'$ (possibly $k$) such that $\theta_{k',i}=1$ and $\tilde\zeta^Y_{k',j'}\neq 0$ is part of the slow process.

\end{condition}

This condition is simplest to explain in case the jump rates on the time scale of the fast dynamics  are given by mass-action functions, that is, a reaction $k$ acting on fast species $\{S_i\}$:
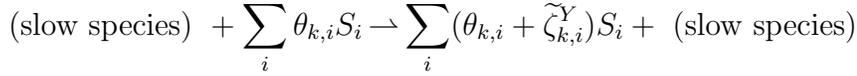
\[ \mbox{(slow species) }+\sum_i\theta_{k,i}S_i \mathop{\rightharpoonup}\sum_i (\theta_{k,i}+\tilde\zeta^Y_{k,i})S_i + \mbox{ (slow species) }\]
has jump rate of the form:
\[\lambda_k(z)=\lambda_k(x)\mathop{\prod}\limits _i y^{ \theta_{k,i}}.\]
Linearity of the jump rates (in the fast variables) is implied by $\theta_{k,i}\in \{0,1\}, \forall k,i$. 
When a reaction has its jump's change for species $i$ is  $\tilde\zeta_{k,i}^Y=-1$ then this is a reaction which is either simply using up species $i$ or converting species $i$ into some other species $j\neq i$: 
\[ \mbox{(slow species) }+S_i \mathop{\rightharpoonup} \mbox{ (slow species)} + \tilde\zeta^Y_{k,j}S_j\]
with $\tilde\zeta^Y_{k,j}=0$ in the former and $\tilde\zeta^Y_{k,j}>0$ in the latter case. Since species $i$ goes into this reaction, it implies that also $\theta_{k,i}=1$. There may be other reactions $k'$ with $\theta_{k',i}=1$ which only use species $i$  as a catalyst and then $\tilde\zeta^Y_{k',i}\ge 0$:
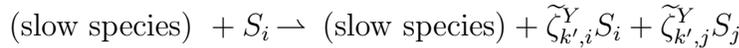
\[ \mbox{(slow species) }+S_i \mathop{\rightharpoonup} \mbox{ (slow species)} + \tilde\zeta^Y_{k',i}S_i+ \tilde\zeta^Y_{k',j}S_j\]
Our condition says that any species $j\neq i$ created by such a reaction from $i$ must be either in a conservation relation with it: $\sum_i\tilde\theta_iY^N_i(t)=\sum_i\tilde\theta_iY^N_i(0)$ for all $t>0$, or is a slow species. 
For simplicity we assume that if species $i$  is part of a conservation law $\tilde\theta$ with $\tilde\theta_i>0$ then there is a unique species $i'$ such that $\tilde\theta_{i'}>0$ as well, we denote this relationship by $i'\sim i$. Let $M_i=\tilde\theta_iY^N_i(t)+\tilde\theta_{i'}Y^N_{i'}(t), \forall t\ge 0$, and in such pairs apply a change of variables  $y_{i'}=M_i/\tilde\theta_{i'}-(\tilde\theta_{i}/\tilde\theta_{i'})y_i$.



We can now solve the eigenvalue problem when the fast process is purely a Markov chain.
\begin{proposition}\label{prop-EVP}
Assume the effective dynamics of the fast process is a Markov chain and Condition~\ref{cond-conversion} holds. Let ${\mathcal I}=\{i: \exists k,j\; \theta_{k,i}= 1,\tilde\zeta^Y_{k,j}=-1\}, \mathcal I_j=\{i:\exists k\; \theta_{k,i}=1, \tilde\zeta^Y_{k,j}=-1\}$, and ${\mathcal J}=\{j: \exists k\;\tilde\zeta^Y_{k,j}=-1\}$. Let $\{A_{ij},B_i,C_{ij}\}$ be the functions given in \eqref{eq-Ax}-\eqref{eq-Cx}. 
If $\forall x\in E, p\in\mathbb R$:\; {\bf (a)} for each $j\in\mathcal J$ the quadratic equation
\[z_j^2\sum_{i\in\mathcal I_j}A_{ij}(x,p)+z_j\sum_{i\in\mathcal I_j}B_{ij}(x,p)+\sum_{i\in\mathcal I_j}C_{ij}(x,p)=0\] 
has a unique positive solution $z_j$; and:\; {\bf (b)} over the set ${\mathcal I}^c=\{i:  \nexists k,j\; \theta_{k,i}= 1,\tilde\zeta^Y_{k,j}=-1\}$  the system of linear equations
\[\Big\{\sum_{j\in\mathcal J^c} z_jA_{ij}(x,p)+\sum_{j\in\mathcal J^c}B_{ij}(x,p)=0\Big\}_{i\in {\mathcal I}^c}\]
has a unique positive solution of variables $\{z_j\}_{j\in\mathcal J^c}$ over the set $\mathcal J^c=\{j: \nexists k\;\tilde\zeta^Y_{k,j}=-1\}$.
Then, for any $j\in\mathcal J$ for which $\sum_{i\in\mathcal I_j}A_{ij}(x,p)\neq 0$ only {\bf (i)} in Condition~\ref{cond-conversion} is possible. Moreover, 
the eigenvalue problem \eqref{eq-EVP} has a unique solution given by 
\begin{eqnarray}\label{eq-EVPsoln}
\ol H_0(x,p)&=&\!\!\!\sum_{k:\theta_{k,0}=1}\lambda_k(x)\Big({\mathbf 1}_{\beta_k=1}(e^{\tilde \zeta_k^X\cdot p}e^{\sum_i a_i(x)\tilde \zeta^Y_{k,i}}-1)+{\mathbf 1}_{\beta_k>1}\tilde \zeta^X_k\cdot p\Big)\\
&+&\sum_i {\mathbf 1}_{\exists i' \sim i}(M_i/\tilde\theta_{i'})\sum_{k:\theta_{k,i}=1}\lambda_k(x)\Big({\mathbf 1}_{\beta_k=1}(e^{\tilde \zeta_k^X\cdot p}e^{\sum_i a_i(x)\tilde \zeta^Y_{k,i}}-1)+{\mathbf 1}_{\beta_k>1}\tilde \zeta^X_k\cdot p\Big)\nonumber
\end{eqnarray}
where $e^{\sum_i a_i(x)\tilde \zeta^Y_{k,i}}=\prod_i z_i^{\tilde \zeta^Y_{k,i}}$ in the formula is determined  by the solutions to the quadratic and linear equations above, and the associated eigenfunction is given by $e^{h(x,y)}=\prod_i z_i^{y_i}$.
\end{proposition}

\begin{proof} We rewrite (\ref{eq-EVP}) with $h(x,y)=\sum_ia_i(x)y_i$
\begin{eqnarray*}
&&\!\!\!\!\!\!\!\!\!\!\left(V(y;x,p)+e^{-h(x,y)} L_1^{x,p} e^{h(x,y)}\right)\\
=&& \sum_{k:\beta_k=1}\lambda_k(z)(e^{\tilde \zeta_k^X\cdot p+\sum_i a_i(x)\tilde \zeta^Y_{k,i}}-1)+\sum_{k:\beta_k>1}\lambda_k(z)\tilde \zeta^X_k\cdot p\\
=&& \sum_{k:\beta_k=1}\lambda_k(x)(\theta_{k,0}+\sum_i\theta_{k,i} y_i)(e^{\tilde \zeta_k^X\cdot p+\sum_i a_i(x)\tilde \zeta^Y_{k,i}}-1)+\sum_{k:\beta_k>1}\lambda_k(x)(\theta_{k,0}+\sum_i\theta_{k,i} y_i)\tilde \zeta^X_k\cdot p\\
=&&\sum_{k:\theta_{k,0}=1}\lambda_k(x)\Big({\mathbf 1}_{\beta_k=1}(e^{\tilde \zeta_k^X\cdot p+\sum_i a_i(x)\tilde \zeta^Y_{k,i}}-1)+{\mathbf 1}_{\beta_k>1}\tilde \zeta^X_k\cdot p\Big)\\
&&+\!\!\sum_{i\in\mathcal I^c}
y_i\sum_{k:\theta_{k,i}=1}\lambda_k(x)\Big({\mathbf 1}_{\beta_k=1}(e^{\tilde \zeta_k^X\cdot p+\sum_i a_i(x)\tilde \zeta^Y_{k,i}}-1)+{\mathbf 1}_{\beta_k>1}\tilde \zeta^X_k\cdot p\Big)\\
&&+\!\!\sum_{i\in\mathcal I}
y_i\sum_{k:\theta_{k,i}=1}\lambda_k(x)\Big({\mathbf 1}_{\beta_k=1}(e^{\tilde \zeta_k^X\cdot p+\sum_i a_i(x)\tilde \zeta^Y_{k,i}}-1)+{\mathbf 1}_{\beta_k>1}\tilde \zeta^X_k\cdot p\Big)
\end{eqnarray*}

In order to get a result which is independent of $y$ we need to set all of the terms next to $y_i$ to zero, which will leave the first row intact, except for adding to it all the $M_i$ dependent terms from the change of variables  $y_{i'}=M_i/\tilde\theta_{i'}-(\tilde\theta_{i}/\tilde\theta_{i'})y_i$ where $i'\sim i$ exists 
\begin{eqnarray*}
\sum_i {\mathbf 1}_{\exists i' \sim i}(M_i/\tilde\theta_{i'})\sum_{k:\theta_{k,i}=1}\lambda_k(x)\Big({\mathbf 1}_{\beta_k=1}(e^{\tilde \zeta_k^X\cdot p+\sum_i a_i(x)\tilde \zeta^Y_{k,i}}-1)+{\mathbf 1}_{\beta_k>1}\tilde \zeta^X_k\cdot p\Big)
\end{eqnarray*}

Using  Condition~\ref{cond-conversion} we can rewrite the middle row as
\begin{eqnarray*}
\sum_{i\in\mathcal I^c}
y_i\sum_{j\in\mathcal J^c}\big(e^{a_j(x)}A_{ij}(x,p)+B_{ij}(x,p)\big),
\end{eqnarray*}
 and likewise the last  row as 
\begin{eqnarray*}
\sum_{i\in\mathcal I}
y_i\sum_{j\in\mathcal J}\big(e^{a_j(x)}A_{ij}(x,p)+B_{ij}(x,p)+e^{-a_j(x)}C_{ij}(x,p)\big).
\end{eqnarray*}
where we used the following functions 
\begin{eqnarray}
A_{ij}(x,p)&=&\!\!\!\sum_{k:\theta_{k,i}=1}{\mathbf 1}_{\tilde\zeta^Y_{k,j}=1}\lambda_k(x){\mathbf 1}_{\beta_k=1}e^{\tilde\zeta^X_k\cdot p}-{\mathbf 1}_{i' \sim i}({\tilde\theta_i}/{\tilde\theta_{i'}})\sum_{k:\theta_{k,i'}=1}{\mathbf 1}_{\tilde\zeta^Y_{k,j}=1}\lambda_k(x){\mathbf 1}_{\beta_k=1}e^{\tilde\zeta^X_k\cdot p}\label{eq-Ax}\\
B_{ij}(x,p)&=&\!\!\!\sum_{k:\theta_{k,i}=1}{\mathbf 1}_{\tilde\zeta^Y_{k,j}=0}\lambda_k(x){\mathbf 1}_{\beta_k=1}e^{\tilde\zeta^X_k\cdot p}-{\mathbf 1}_{i' \sim i}({\tilde\theta_i}/{\tilde\theta_{i'}})\sum_{k:\theta_{k,i'}=1}{\mathbf 1}_{\tilde\zeta^Y_{k,j}=0}\lambda_k(x){\mathbf 1}_{\beta_k=1}e^{\tilde\zeta^X_k\cdot p}\label{eq-Bx}\\
&+&  \!\!\!\sum\limits_{k:\theta_{k,i}=1}\lambda_k(x)\big(-{\mathbf 1}_{\beta_k=1}+{\mathbf 1}_{\beta_k>1}\tilde\zeta^X_k\!\!\cdot p\big)-{\mathbf 1}_{i' \sim i}({\tilde\theta_i}/{\tilde\theta_{i'}})\sum_{\theta_{k,i'}=1}\lambda_k(x)\big(-{\mathbf 1}_{\beta_k=1}+{\mathbf 1}_{\beta_k>1}\tilde\zeta^X_k\!\!\cdot p\big)\nonumber\\
C_{ij}(x,p)&=&\!\!\!\sum_{k:\theta_{k,i}=1}{\mathbf 1}_{\tilde\zeta^Y_{k,j}=-1}\lambda_k(x){\mathbf 1}_{\beta_k=1}e^{\tilde\zeta^X_k\cdot p}-{\mathbf 1}_{i' \sim i}({\tilde\theta_i}/{\tilde\theta_{i'}})\sum_{k:\theta_{k,i'}=1}{\mathbf 1}_{\tilde\zeta^Y_{k,j}=-1}\lambda_k(x){\mathbf 1}_{\beta_k=1}e^{\tilde\zeta^X_k\cdot p}\label{eq-Cx}
\end{eqnarray}
For each species $i$ from the middle row we get one equation in a system of $|\mathcal I^c |$ linear equations in the variables $z_j=e^{a_j(x)}$ over the set $j\in \mathcal J^c$. 
For each $i$ from the last row we get a single quadratic equation in the variables $z_j=e^{a_j(x)}$, note that Condition~\ref{cond-conversion} implies the map $i\mapsto j$ is unique. The set $\mathcal I_j=\{i:\exists k\; \theta_{k,i}=1, \tilde\zeta_{k,j}=-1\}$ may not be of size $1$, so these quadratic equations will combine to identify a single solution $z_j$ to the quadratic equation
\[\sum_{i\in \mathcal I_j}A_{ij}(x)e^{a_j(x)}+\sum_{i\in \mathcal I_j}B_{ij}(x)+\sum_{i\in \mathcal I_j}C_{ij}(x)e^{-a_j(x)}=0.\]
Note that this equation can in fact have the coefficient $\sum_{i\in \mathcal I_j}A_{ij}=0$ in which case option {\bf (ii)} of the Condition~\ref{cond-conversion} is possible.
However, if the equation has $\sum_{i\in \mathcal I_j}A_{ij}\neq 0$ then the quadratic equation has a unique positive solution iff $\sum_{i\in\mathcal I_j}A_i(x)\sum_{i\in\mathcal I_j}C_i(x)<0$. This implies that only option {\bf (i)} of the Condition~\ref{cond-conversion} is possible and species $j\in\mathcal J$ has to be  in a conservation law with some other species, as otherwise all of the functions $A_{ij}(x,p), C_{ij}(x,p)$ are positive and only $B_{ij}(x,p)$ can be negative.  
\end{proof}

Note that $|\mathcal I|\ge |\mathcal J|$ so $|\mathcal J^c|\ge |\mathcal I^c|$ and in general there may be variables that are not defined by the system (e.g. $\emptyset \rightharpoonup S$ with no $S\rightharpoonup$ in system). This is the reason we include in out state space $Z$ only {\it active} species in the reaction network which are defined as species that appears as an input in at least one reaction. An example of an inactive species is the  product species $P$ in the enzymatic kinetics example, or the packaged virus particle species $V$ in the viral production example.\\

We next provide an analogous result in case the effective dynamics of the fast process is a piecewise-deterministic  Markov chain and a dynamical system (a PDMP). In this case sufficient conditions for the existence of an explicit solution to the eigenvalue problem is slightly more messy to state, but equally straightforward to solve. Recall that the state space of a PDMP separates into a discrete and a continuous component, with the discrete component jumps according to Markov chain and the continuous component performs deterministic dynamics. In case the fast process is a PDMP, no restrictions are necessary on its continuous component, and conditions on its discrete component are similar to the case when the fast process is just a Markov chain. However, here we have an extra option (iii) in which species created from reactions using fast discrete species are allowed to be fast continuos species (similarly to the option (ii) where they are allowed to be slow species, see the explanation after Condition~\ref{cond-conversion}). 

\begin{condition}\label{cond-conversionPDMP}
Suppose there are constants $\theta_{k,0},\{\theta_{k,i}\}$ $\in \{0,1\}$ satisfying $\theta_{k,0}+\sum\theta_{k,i}=1$ such that
 \[\lambda_k(x,y) = \lambda_k(x)(\theta_{k,0}+\sum_i\theta_{k,i} y_i);\]
for each  $i$ in the {\bf discrete} component $\tilde\zeta^Y_{k,i} \in\{-1,0,1\}\; \forall k$; and for each  $j$ in the {\bf discrete} component for which $\exists k,i$ such that $\theta_{k,i}= 1$ and $\tilde\zeta^Y_{k,j}=-1$ we have that either 
{\bf (i)} 
or 
{\bf (ii)} 
of Condition~\ref{cond-conversion} hold 
or:

 {\bf (iii)} a species $j'\neq j$  for which $\exists k'$ (possibly $k$) such that $\theta_{k',i}=1$ and $\tilde\zeta^Y_{k',j'}\neq 0$ is in the continuous component of the fast process, and is such that: $\exists i^\ast\neq i, \exists k^\ast\neq k'$ such that $\theta_{k^\ast,i^\ast}=1$  and $\tilde\zeta^Y_{k^\ast,j'}\neq 0$  and also $\tilde\zeta^Y_{k^\ast,j^\ast}\neq 0$ only if $j^\ast$ is a slow species.
\end{condition}

We note that if {\bf (i)} holds for species $i$ than the conserved species $i'\sim i$ must also be in the discrete component of the fast process. 

\begin{corollary}\label{cor-EVPPDMP} 
Assume the effective dynamics of the fast process is a piecewise-deterministic Markov chain and Condition~\ref{cond-conversionPDMP} holds. Let $\mathcal I,\mathcal I_j, \mathcal J$ be as in Proposition~\ref{prop-EVP}. If $\forall x\in E, p\in\mathbb R$ for each discrete component species $j\in Y_d\cap\mathcal J$  the quadratic equation
\[z_j^2\sum_{i\in\mathcal I_j}A_{jj}(x,p)+z_j\sum_{i\in\mathcal I_j}B_{ij}(x,p)+\sum_{i\in\mathcal I_j}C_{ij}(x,p)=0\] 
has a unique positive solution $z_j$, and over the set $\mathcal I^c$ the system of linear equations 
\[\Big\{\sum_{j\in Y_d\cap\mathcal J^c} z_jA_{ij}(x,p)+\sum_{j\in Y_c} u_j A^c_{ij}(x,p)+\sum_{j\in (Y_d\cap\mathcal J^c)\cup Y_c}B_{ij}(x,p)\}=0\Big\}_{i\in \mathcal I^c}\]
has a unique positive solution of variables $\{z_j\}_{j\in Y_d\cap\mathcal J^c}, \{u_j\}_{j\in Y_c}$, 
then for any $j\in Y_d\cap \mathcal J$ for which $\sum_{i\in\mathcal I_j}A_{ij}(x,p)\neq 0$ only {\bf (i)} in Condition~\ref{cond-conversionPDMP} is possible. 
Moreover, the eigenvalue problem \eqref{eq-EVP} has a unique solution given by 
\begin{eqnarray}\label{eq-EVPsolnPDMP}
\ol H_0(x,p)&=&\!\!\!\sum_{k:\theta_{k,0}=1}\lambda_k(x)\Big({\mathbf 1}_{\beta_k=1}(e^{\tilde \zeta_k^X\cdot p}e^{\sum_i a_i(x)\tilde \zeta^Y_{k,i}}-1)+{\mathbf 1}_{\beta_k>1}(\tilde \zeta^X_k\cdot p+\tilde \zeta^Y_k\cdot a(x))\Big)\\
&+&\sum_i {\mathbf 1}_{\exists i' \sim i}(M_i/\tilde\theta_{i'})\sum_{k:\theta_{k,i}=1}\lambda_k(x)\Big({\mathbf 1}_{\beta_k=1}(e^{\tilde \zeta_k^X\cdot p}e^{\sum_i a_i(x)\tilde \zeta^Y_{k,i}}-1)+{\mathbf 1}_{\beta_k>1}(\tilde \zeta^X_k\cdot p+\tilde \zeta^Y_k\cdot a(x))\Big)\nonumber
\end{eqnarray}
where in the above formula the terms $e^{\sum_i a_i(x)\tilde \zeta^Y_{k,i}}=\prod_{i\in Y_d}z_i^{\tilde \zeta^Y_{k,j}}$ and  $\tilde \zeta^Y_k\cdot a(x)=\sum_{j\in Y_c}\tilde\zeta^Y_{k,j}u_j(x)$ are determined  by the solutions to the stated quadratic and linear equations, and the associated eigenfunction is given by $e^{h(x,y)}=\prod_{i\in Y_d} z_i^{y_i}\prod_{i\in Y_c}e^{u_j(x)y_j}$.
\end{corollary}

\begin{proof} 
The generator of the fast process has an additional term from the continuous component
\[L_1^{x,p}f(z)=\sum_{k:\beta_k=1} \lambda_k(z)(f(x,y+\tilde\zeta^Y_k)-f(x,y))+ \sum_{k:\beta_k>1} \lambda_k(z) \tilde\zeta^Y_k \cdot\nabla_Y f(x,y)\]
so if $a(x)=\{a_i(x)\}$ using $h(x,y)=\sum_ia_i(x)y_i$ as before the equation \eqref{eq-EVP} now becomes
\begin{eqnarray*}
&&\!\!\!\!\!\!\!\!\!\!\left(V(y;x,p)+e^{-h(x,y)} L_1^{x,p} e^{h(x,y)}\right)\\
=&&\sum_{k:\theta_{k,0}=1}\lambda_k(x)\Big({\mathbf 1}_{\beta_k=1}(e^{\tilde \zeta_k^X\cdot p+\sum_i a_i(x)\tilde \zeta^Y_{k,i}}-1)+{\mathbf 1}_{\beta_k>1}(\tilde \zeta^X_k\cdot p+\tilde \zeta^Y_k\cdot a(x))\Big)\\
&&+\!\!\sum_{i\in\mathcal I^c}
y_i\sum_{k:\theta_{k,i}=1}\lambda_k(x)\Big({\mathbf 1}_{\beta_k=1}(e^{\tilde \zeta_k^X\cdot p+\sum_i a_i(x)\tilde \zeta^Y_{k,i}}-1)+{\mathbf 1}_{\beta_k>1}(\tilde \zeta^X_k\cdot p+\tilde \zeta^Y_k\cdot a(x))\Big)\\
&&+\!\!\sum_{i\in\mathcal I}
y_i\sum_{k:\theta_{k,i}=1}\lambda_k(x)\Big({\mathbf 1}_{\beta_k=1}(e^{\tilde \zeta_k^X\cdot p+\sum_i a_i(x)\tilde \zeta^Y_{k,i}}-1)+{\mathbf 1}_{\beta_k>1}(\tilde \zeta^X_k\cdot p+\tilde \zeta^Y_k\cdot a(x))\Big)
\end{eqnarray*}
where we rewrite the middle row, using $Y_d$ and $Y_c$ to denote discrete and continuous components respectively, as
\begin{eqnarray*}
\sum_{i\in\mathcal I^c}
y_i\big(\sum_{j\in Y_d\cap\mathcal J^c}e^{a_j(x)}A_{ij}(x,p)+\sum_{j\in Y_c}{a_j(x)}A^c_{ij}(x,p)+\sum_{j\in (Y_d\cap\mathcal J^c)\cup Y_c}B_{ij}(x,p)\big),
\end{eqnarray*}
 and likewise the last  row as 
\begin{eqnarray*}
\sum_{i\in\mathcal I}
y_i\sum_{j\in Y_d\cap\mathcal J}\big(e^{a_j(x)}A_{jj}(x,p)+B_{ij}(x,p)+e^{-a_j(x)}C_{ij}(x,p)\big).
\end{eqnarray*}
where the functions $A_{ij}(x),B_i(x),C_i(x)$ are the same as in \eqref{eq-Ax},\eqref{eq-Bx},\eqref{eq-Cx}, and the only new contribution is from the function
\begin{eqnarray}
A^c_{ij}(x,p)&=&\!\!\!\sum_{k:\theta_{k,i}=1}\lambda_k(x){\mathbf 1}_{\beta_k>1}\tilde\zeta^Y_{k,j}-{\mathbf 1}_{i' \sim i}({\tilde\theta_i}/{\tilde\theta_{i'}})\sum_{k:\theta_{k,i'}=1}\lambda_k(x){\mathbf 1}_{\beta_k>1}\tilde\zeta^Y_{k,j}.\label{eq-AcxPDMP}
\end{eqnarray}
Otherwise the rest of the process of solving for the result is the same.
\end{proof}

The condition precludes systems with reactions that are bi-molecular in fast species, which is an assumption we made for the sake of simplicity of the proposition. It is in principle clear how one can try to extend the above derivation in case of multi-molecular reactions between fast species to get  a solution to the eigenvalue problem. 

We will show how the procedure for solving the eigenvalue problem works in a variety of examples, including the two mentioned in the Introduction.\\

\def\ep{\epsilon}
\def\N{{\mathbb{N}}}
\def\R{{\mathbb{R}}}
\def\Z{{\mathbb{Z}}}
\def\l{\left}
\def\r{\right}
\def\p{\partial}
\def\L{\mathcal L}
\def\ol{\overline}
\def\ul{\underline}
\def\ta{\tilde a}

\renewcommand {\theequation}{\arabic{section}.\arabic{equation}}
\def \non{{\nonumber}}
\def \hat{\widehat}
\def \tilde{\widetilde}
\def \bar{\overline}

\setcounter{equation}{0}

\section{Examples}\label{sec-examples}

We assume that reaction rates have mass-action form throughout the following examples, as in \eqref{eq-massact} although in one of the examples (Subsection~\ref{subsec-SRG}) we will in addition allow some of the chemical reaction constants $\kappa_i$ to be a function of a specific species of interest in the system. We consider the LDP on $E_X=(\ep,\infty)$  for arbitrarily small $\ep>0$, and $E_Y$ is either $[0,K]$ for some $K$ or $(0,\infty)$ depending on existence of conservation laws in each example. 

\subsection{Enzymatic kinetics (MM)}\label{subsec-MM}

We recall the model for enzyme kinetics (Michealis-Menten), with an inflow of the substrate:
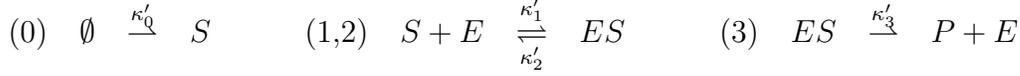
\begin{center}
\begin{tabular}{lrlllrlllrlll} 
(0)&$\emptyset$&$\stackrel {\kappa'_{0}}{\rightharpoonup}$& $S$&\qquad
(1,2)&$S+E$&$\mathop{\stackrel {\kappa'_{1}}{\rightleftharpoons}}\limits_{\kappa'_2}$&$ES$&\qquad
(3)&$ES$&$\stackrel {\kappa'_{3}}{\rightharpoonup}$&$P+E$\\
\end{tabular}
\end{center}
The scaling of the amounts is implied by the fact that molecular amount of the substrate $S$ is an order of magnitude greater than the amount of enzyme $E$ and of the enzyme-subtsrate complex $ES$.
Let $Z_1, Z_2, Z_3, Z_4$ represent the amounts of $S,E,ES,P$ molecular species respectively, the orders of magnitude lead to appropriately scaled species amounts:  $Z^N_1=Z_1/N, Z^N_2=Z_2,Z^N_3=Z_3, Z^N_4=Z_4/N$.
The reaction constants $\kappa_i$ also have different orders of magnitude, with those for the dissolution of the enzyme-substrate complex in reactions (2,3) being an order of magnitude larger than the forming of the complex. Let $\kappa_0=\kappa'_0/N,\kappa_1=\kappa'_1, \kappa_2=\kappa'_2/N, \kappa_3=\kappa'_3/N$, so the  model of the system is 
\begin{eqnarray*}
Z^N_1(t)&=&Z^N_1(0)+-N^{-1}Y_1(N\kappa_0t)-N^{-1}Y_1(N\int_0^t\kappa_1Z^N_1(s)Z^N_2(s)ds
)\\&&\qquad\quad+N^{-1}Y_2(N\int_0^t\kappa_2(M-Z_2^N(s))ds)\\
Z^N_2(t)&=&Z^N_2(0)-Y_1(N\int_0^t\kappa_1Z^N_1(s)Z^N_2(s)ds)+Y_2(
N\int_0^t\kappa_2(M-Z_2^N(s))ds) \\&& \qquad\quad+Y_3(N\int_0^t\kappa_3(M-Z_2^N(s))ds)\\
Z_3^N(t)&=&Z^N_3(0)+Y_1(N\int_0^t\kappa_1Z^N_1(s)Z^N_2(s)ds)-Y_2(N\int_0^t\kappa_2Z^N_3(s)ds) -Y_3(N\int_0^t\kappa_3Z^N_3(s)ds)\\
Z^N_4(t)&=&N^{-1}Y_3(N\int_0^t\kappa_3Z^N_3(s)ds),
\end{eqnarray*}
There is a conservation law between $E$ and $ES$ since $Z_2^N(t)+Z_3^N(t)\equiv M, \forall t>0$ hence we will use a change of variables $z_3=M-z_2$. We also note that $P$ is not an ``active" species of the system, as it doesn't enter on the left hand side of any reaction.

This leaves a system with the slow and fast process: $X^N=Z_1^N$ and $Y^N=Z_2^N$ respectively. The scaling conditions are clearly satisfied with the time-scale separation between slow $X^N$ and $Y^N$ of order $N$. The effective dynamics of the slow process is given by the ODE in \eqref{eq-MMlim} while the effective dynamics of the fast process is a birth-death Markov chain with death rate $\kappa_1z_1z_2$ and birth rate $(\kappa_2+\kappa_3)(M-z_2)$.

The rate of reaction (1) is binary, with $\lambda_1(z)=\kappa_1z_1z_2$, but since one of the factors, $Z_2$, is part of a conservation law and hence bounded by the constant $M$, by the Proposition~\ref{prop-binbound} the process has an exponentially good approximation in a sequence of processes, indexed by increasing values of $M'$, in which the rate $\lambda_1$ is replaced by:
$\lambda'_1(z)=\kappa_1z_1(z_2\wedge M')\equiv \lambda_1(z)$ for $M'\ge M$. 

The generator for the pair is:
\begin{eqnarray*}
A_Nf(x,y)&=&N\kappa_0(f(x+N^{-1},y)-f(x,y))\\&+&N\kappa_1xy(f(x-N^{-1},y-1)-f(x,y))\\&+&N\kappa_2(M-y)(f(x+N^{-1},y+1)-f(x,y))\\
&+&N\kappa_3(M-y)(f(x,y+N^{-1})-f(x,y))
\end{eqnarray*}
so the exponential generator, acting on $f_N(x,y)=f(x)+N^{-1}h(x,y)$, and its limit
are:
\begin{eqnarray*}
H_Nf_N(x,y)&=&\kappa_0(e^{N(f(x+N^{-1})-f(x))+h(x+N^{-1},y)-h(x,y)}-1) \\
&+&\kappa_1xy(e^{N(f(x-N^{-1})-f(x))+h(x-N^{-1},y-1)-h(x,y)}-1) \\&+&\kappa_2(M-y)(e^{N(f(x+N^{-1})-f(x))+h(x+N^{-1},y+1)-h(x,y)}-1)\\
&+&\kappa_3(M-y)(e^{h(x,y+1)-h(x,y)}-1)
\end{eqnarray*}
\begin{eqnarray*}
\lim_{N\rightarrow\infty}H_Nf_N(x,y)
&=&\kappa_0e^{f'(x)}+\kappa_1xy(e^{-f'(x)}e^{h(x,y-1)-h(x,y)}-1)+(\kappa_2+\kappa_3)(M-y)(e^{f'(x)}e^{h(x,y+1)-h(x,y)}-1)\\
&=&\ol H_0(x,f'(x))
\end{eqnarray*}
so that $\ol H_0(x,p)$ satisfies the EVP equation \eqref{eq-EVP} with 
\begin{eqnarray*}
V(y;x,p)&=&\kappa_0e^{p}+\kappa_1xy(e^{-p}-1)+\kappa_2(M-y)(e^{p}-1)\\
L_1^{x,p}e^{g(x,y)}&=&\kappa_1xye^{-p}(e^{g(x,y-1)}-e^{g(x,y)})+(\kappa_2+\kappa_3)(M-y)(e^{g(x,y+1)}-e^{g(x,y)})
\end{eqnarray*}

Since the perturbed effective fast dynamics is a birth-death Markov chain with birth rate $(\kappa_2+\kappa_3)(M-y)$ and death rate $\kappa_1xye^{-p}$ the density Condition~\ref{cond-ergod} is satisfied, and there is a unique stationary distribution which is a Binomial$(M,\pi^{x,p})$ with size parameter $M$ and probability parameter $\pi^{x,p}=(\kappa_2+\kappa_3)/(\kappa_2+\kappa_3+\kappa_1xe^{-p})$.
Since the state space for $Y^{x,p}$ is $\{0,1,\dots,M\}$ both Lyapunov Conditions~\ref{cond-multiplyapun} and ~\ref{cond-multiplyapun} are trivially satisfied.

In order to solve the EVP for $\ol H_0(x,p)$ above let $g(x,y)=a(x)y$, then:
\begin{eqnarray*}
\ol H_0(x,y)&=& \kappa_1xy(e^{-p}e^{-a(x)}-1)+(\kappa_2+\kappa_3)(M-y)(e^{a(x)}-1)\\
\!\!\!&=&y\left(\kappa_1x(e^{-p}e^{-a(x)}-1)-\kappa_2(e^{p}e^{a(x)}-1)-\kappa_3(e^{a(x)}-1)\right)\\&&\quad\quad+M\left(\kappa_2(e^{p}e^{a(x)}-1)+\kappa_3(e^{a(x)}-1)\right)+\kappa_0(e^p-1)
\end{eqnarray*}
which after setting the coefficient of $y$  to zero gives a single quadratic equation:
\[(\kappa_2e^{p}+\kappa_3 )e^{2a(x)}+(\kappa_1x-\kappa_2-\kappa_3)e^{a(x)}-\kappa_1xe^{-p}=0\]
Since there was a conservation law in the fast variable the quadratic coefficients can satisfy $A(x)C(x)<0$ and produce a unique positive solution
\[e^{a(x)}=\frac{-(\kappa_1x-\kappa_2-\kappa_3)+\sqrt{(\kappa_1x-\kappa_2-\kappa_3)^2+4(\kappa_2e^{p}+\kappa_3 )\kappa_1xe^{-p}}}{2(\kappa_2e^{p}+\kappa_3 )}\]
which by Proposition~\ref{prop-EVP} implies
\[\ol H_0(x,p)=\frac M2\left(-\kappa_2-\kappa_3-\kappa_1x+\sqrt{(\kappa_2+\kappa_3-\kappa_1x)^2+4(\kappa_2+\kappa_3e^{-p} )\kappa_1x}\right)+\kappa_0(e^p-1)\]
From the explicit equation it is clear that on $(\ep,\infty)\times \mathbb R$ the function $\ol H_0(x,p)$ is convex and coercive in $p$ and it satisfies the condition \eqref{eq:H2ploc}.\\

\subsection{Self-regulated gene expression (SRG)}\label{subsec-SRG}
Another common example from systems biology is that of self-regulated gene expression

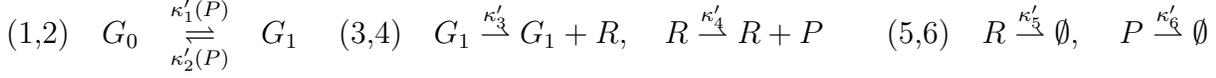
\begin{center}
\begin{tabular}{lrlllrlllrlll} 
(1,2)&$G_0$&$\mathop{\stackrel {\kappa'_1(P)}{\rightleftharpoons}}\limits_{\kappa'_2(P)}$&$G_1$&
\;(3,4)&$G_1\stackrel {\kappa'_{3}}{\rightharpoonup}G_1+R$,\,& $R\stackrel {\kappa'_{4}}{\rightharpoonup}R+P$&\quad (5,6)&$R\stackrel {\kappa'_{5}}{\rightharpoonup}\emptyset$,\; &$P\stackrel {\kappa'_{6}}{\rightharpoonup}\emptyset$\\
\end{tabular}
\end{center}

\noindent where $G_0,G_1$ are inactive and active molecular conformations of a gene, $R$ is the MRNA, and $P$ is a protein expressed by this gene. The self-regulation of the gene comes via the protein it expresses, whose amount affects the rate at which the active conformation becomes an inactive conformation and vice versa.
 
We assume there is one copy of the gene, and the amount of MRNA and protein in the long run is expressed in terms of a scaling parameter $N$. The two conformations of the gene are rapidly changing their states, and the active one is rapidly involved in expression,  so the rates $\kappa_1',\kappa_2', \kappa_3'$ are assumed to be of order $N$, while the translation rate and the MRNA and protein degradation rates $\kappa'_4,\kappa_5', \kappa'_6$ are of order $1$. 
Let $Z_1,Z_2, Z_3, Z_4$ represent the amounts of $G_0,G_1,R,P$ respectively, so their rescaled versions are $Z^N_1=Z_1,Z^N_2=Z_2,Z^N_3=Z_3/N, Z^N_4=Z_4/N$. Let $\kappa_1(Z_3/N)=\kappa_1'(Z_3)/N,\kappa_2(Z_3/N)=\kappa_2'(Z_3)/N,\kappa_3=\kappa_3'/N$ and $\kappa_4=\kappa_4', \kappa_5=\kappa_5', \kappa_6=\kappa_6'$. We assume the self-regulating rates $\kappa_1(\cdot), \kappa_2(\cdot)$ are Lipschitz and grow at most linearly, and to prevent absorption of the system we will assume they are positive.

There is a conservation law between $G_0$ ad $G_1$ since $G_0(t)+G_1(t)\equiv 1, \forall t$ so we use $z_2=1-z_1$. The model of the system is
\begin{eqnarray*}
Z_1^N(t)&=&Z_1^N(0)-{Y}_1(\int_0^t N\kappa_1(Z_4^N(s))Z_1^N(t)ds)+{Y}_2(\int_0^t N\kappa_2(Z_4^N(s))(1-Z_1^N(s))ds)\\
Z_3^N(t)&=&
N^{-1}{Y}_3(\int_0^t N\kappa_3(1-Z_1^N(s))ds)-N^{-1}{Y}_5(\int_0^t N\kappa_5Z_3^N(s)ds)\\
Z_4^N(t)&=&N^{-1}{Y}_4(\int_0^t N\kappa_4Z_3^N(s)ds)-N^{-1}{Y}_6(\int_0^t N\kappa_6Z_4^N(s)ds)
\end{eqnarray*}

The slow and fast process are $(X_1^N,X_2^N)=(Z_3^N,Z_4^N)$ and $Y^N=Z_1^N$ respectively, with the time-scale separation of order $N$, and the generator of  the pair is
\begin{eqnarray*}
A_Nf(x,y)&=&N\kappa_1(x_2)y(f(x,y-1)-f(x,y))+N\kappa_2(x_2)(1-y)(f(x,y+1)-f(x,y))\\&+&N\kappa_3(1-y)(f(x_1+N^{-1},x_2,y)-f(x,y))+N\kappa_4x_1(f(x_1,x_2+N^{-1},y)-f(x,y))\\&+&N\kappa_5x_1(f(x_1-N^{-1},x_2,y)-f(x,y))+N\kappa_6x_2(f(x_1,x_2-N^{-1},y)-f(x,y))
\end{eqnarray*}

The effective dynamics of the fast process is a Markov chain on $\{0,1\}$ with \mbox{$1\mapsto 0$} rate $\kappa_1(x_2)y$ and $0\mapsto 1$ rate $\kappa_2(x_2)(1-y)$ and a unique stationary distribution that is Bernoulli with probability $\pi^x=\kappa_2(x_2)/(\kappa_2(x_2)+\kappa_1(x_2))$. This implies the effective dynamics of the slow process is given by the system of ODEs
\begin{eqnarray*}
\dot{x}_1(t)&=&\frac{\kappa_3\kappa_1(x_2(t))}{\kappa_2(x_2(t))+\kappa_1(x_2(t))}-\kappa_5x_1(t)\\
\dot{x}_2(t)&=&\kappa_4x_1(t)-\kappa_6x_2(t)\\\end{eqnarray*}

The exponential generator acting on $f_N(x,y)=f(x)+N^{-1}g(x,y)$, and its limit are:  
\begin{eqnarray*}
H_Nf_N(x,y)&=&\kappa_1(x_2)y(e^{g(x,y-1)-g(x,y)}-1) \\&+&\kappa_2(x_2)(1-y)(e^{g(x,y+1)-g(x,y)}-1)\\
&+&\kappa_3(1-y)(e^{N(f(x_1+N^{-1},x_2)-f(x))+g(x_1+N^{-1},x_2,y)-g(x,y)}-1)\\&+&\kappa_4x_1(e^{N(f(x_1,x_2+N^{-1})-f(x))+g(x_1,x_2+N^{-1},y)-g(x,y)}-1)\\
&+&\kappa_5x_1(e^{N(f(x_1-N^{-1},x_2)-f(x))+g(x_1-N^{-1},x_2,y)-g(x,y)}-1)\\
&+&\kappa_6x_2(e^{N(f(x_1,x_2-N^{-1})-f(x))+g(x_1,x_2-N^{-1},y)-g(x,y)}-1)\\
\lim_{N\rightarrow\infty}H_Nf_N(x,y)&=&\kappa_1(x)y(e^{g(x,y-1)-g(x,y)}-1)+\kappa_2(x)(1-y)(e^{g(x,y+1)-g(x,y)}-1)\\
&+&\kappa_3(1-y)(e^{\partial_{x_1}f(x)}-1)+\kappa_4x_1(e^{\partial_{x_2}f(x)}-1)\\
&+&\kappa_5x_1(e^{-\partial_{x_1}f(x)}-1)+\kappa_6x_2(e^{-\partial_{x_2}f(x)}-1)
\end{eqnarray*}
and we can identify $V$ and $L_1^{x,p}$ as 
\begin{eqnarray*}
V(y;x,p)&=&
\kappa_3(1-y)(e^{p_1}-1)+\kappa_4x_1(e^{p_2}-1)+\kappa_5x_1(e^{-p_1}-1)+\kappa_6x_2(e^{-p_2}-1)\\
L_1^{x,p}e^{g(x,y)}&=&\kappa_1(x_2)y(e^{g(x,y-1)}-e^{g(x,y)})+\kappa_2(x_2)(1-y)(e^{g(x,y+1)}-e^{g(x,y)})
\end{eqnarray*}

Since the perturbed effective fast dynamics is a simple $\{0,1\}$ Markov chain with positive transition rates  (same as in the unperturbed case) the density Condition~\ref{cond-ergod} is satisfied, and both Lyapunov Conditions~\ref{cond-lyapun} and ~\ref{cond-multiplyapun} are trivially satisfied.

Letting $g(x,y)=a(x)y$, implies
\begin{eqnarray*}
\ol H_0(x,p)&=& \kappa_1(x_2)y(e^{-a(x)}-1)+\kappa_2(x_2)(1-y)(e^{a(x)}-1)+\kappa_3(1-y)(e^{p_1}-1)\\
&+&\kappa_4x_1(e^{p_2}-1)+\kappa_5x_1(e^{-p_1}-1)+\kappa_6x_2(e^{-p_2}-1)\\
&=&y\big(\kappa_1(x_2)(e^{-a(x)}-1)-\kappa_2(x_2)(e^{a(x)}-1)-\kappa_3(e^{p_1}-1)\big)\\
&+&\kappa_2(x_2)(e^{a(x)}-1)+\kappa_3(e^p-1)+\kappa_4x_1(e^{p_2}-1)+\kappa_5x_1(e^{-p_1}-1)+\kappa_6x_2(e^{-p_2}-1)
\end{eqnarray*}
which after setting the coefficient of $y$ to $0$ gives the quadratic equation 
\begin{eqnarray*}
-\kappa_2(x_2)e^{2a(x)}+(\kappa_2(x_2)-\kappa_1(x_2)-\kappa_3(e^{p_1}-1))e^{a(x)}+\kappa_1(x_2)=0
\end{eqnarray*}
with $A(x)C(x)<0$ and one positive solution (regardless of the functions $\kappa_1(\cdot),\kappa_2(\cdot)$ or the values of reaction coefficients $\kappa_i, i=3,4,5,6$)
\begin{eqnarray*}
e^{a(x)}=\frac{\kappa_2(x_2)-\kappa_1(x_2)-\kappa_3(e^{p_1}-1)+\sqrt{(\kappa_2(x_2)-\kappa_1(x_2)-\kappa_3(e^{p_1}-1))^2+4\kappa_2(x_2)\kappa_1(x_2)}}{2\kappa_2(x_2)}
\end{eqnarray*}
By Proposition~\ref{prop-EVP} 
\begin{eqnarray*}
\ol H_0(x,p)&=&-\kappa_2(x_2)-\kappa_1(x_2)+\sqrt{(\kappa_2(x_2)-\kappa_1(x_2)-\kappa_3(e^{p_1}-1))^2+4\kappa_2(x_2)\kappa_1(x_2)}\\
&&+\kappa_4x_1(e^{p_2}-1)+\kappa_5x_1(e^{-p_1}-1)+\kappa_6x_2(e^{-p_2}-1)
\end{eqnarray*}
Since the rates $\kappa_1(x), \kappa_2(x)$ are assumed Lipschitz in $x$, with at most linear growth, $\ol H_0(x,p)$ is convex and coercive in $p$ and satisfies the condition \eqref{eq:H2ploc} in  $(\ep,\infty)\times \mathbb R$ for arbitrary $\ep>0$. In \cite{LL17} authors obtain the same large deviation principle for this example (note that the slow and fast variables are labeled differently there).\\

\subsection{Down-regulation (DR)}
We consider a simple model of a negative self-regulation mechanism:
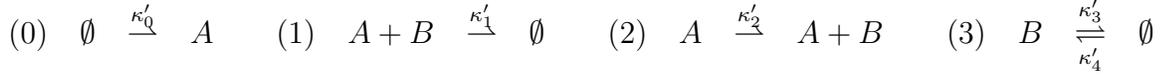
\begin{center}
\begin{tabular}{lrlllrlllrlllrlll} 
(0)&$\emptyset$&$\stackrel {\kappa'_{0}}{\rightharpoonup}$&$A$&
\quad(1)&$A+B$&$\stackrel {\kappa'_{1}}{\rightharpoonup}$&$\emptyset$&
\quad(2)&$A$&$\stackrel {\kappa'_{2}}{\rightharpoonup}$&$A+B$&
\quad(3)&$B$&$\mathop{\stackrel {\kappa'_{3}}{\rightleftharpoons}}\limits_{\kappa'_4}$&$\emptyset$
\\
\end{tabular}
\end{center}
where $A$ is a species of interest and $B$ is a species used to down-regulate it, namely, $A$ controls its own molecular amount by producing more of the regulating species $B$. Suppose the molecular amounts and rates satisfy:
\[|A|=O(N)\mapsto X^N=|A|/N, \; |B|=O(1), \quad\kappa'_{1},\kappa'_{2}\sim O(1),\;  \kappa'_{0},\kappa'_{3}, \kappa'_{4}\sim O(N)\]
Let $Z_1, Z_2$ represent the amounts of $A$ and $B$ molecules respectively, so that the rescaled versions are $Z_1^N=Z_1/N, Z_2^N=Z_2$. Let $\kappa'_{1}=\kappa_1,\kappa'_{2}=\kappa_2$ and $\kappa_3=\kappa'_{3}N, \kappa_4=\kappa'_{4}N$. The model for the system is
\begin{eqnarray*}
Z_1^N(t)&=&Z_1^N(0)+N^{-1}Y_0(\kappa_0t)-N^{-1}Y_1(\int_0^t N\kappa_1 Z_1^N(s)Z_2^N(s)ds)\\
Z_2^N(t)&=&Z_2^N(0)- Y_1(\int_0^t N\kappa_1 Z_1^N(s)Z_2^N(s)ds)
+ Y_2(\int_0^t N\kappa_2 Z_1^N(s)ds)\\&&\quad\qquad -Y_3(\int_0^t N\kappa_3 Z_2^N(s)ds)+ Y_4(\int_0^t N\kappa_4 ds)\\
\end{eqnarray*}
The slow and the fast process are $X^N=Z_1^N$ and $Y^N=Z_2^N$ respectively, and the generator of the  pair is 
\begin{eqnarray*}
A_Nf(x,y)&=&N\kappa_0(f(x+N^{-1},y)-f(x,y))+N\kappa_1yx(f(x-N^{-1},y-1)-f(x,y))\\&+&N\kappa_2x(f(x,y+1)-f(x,y))+N\kappa_3y(f(x,y-1)-f(x,y))\\&+&N\kappa_4(f(x,y+1)-f(x,y))
\end{eqnarray*}
The effective dynamics of the fast process is a simple birth and death Markov chain on $\mathbf Z_+$ with birth rate $\kappa_2x+\kappa_4$ and death rate $\kappa_1xy+\kappa_3y$ with a unique stationary distribution that is Poisson($\mu^x$) with parameter $\mu^x=(\kappa_2x+\kappa_4)/(\kappa_1x+\kappa_3)$. The effective dynamics of the slow process $x$ is a solution to the ODE
\[\dot{x}(t)=\kappa_0-\kappa_1x(t) \frac{\kappa_2x(t)+\kappa_4}{\kappa_1x(t)+\kappa_3}\]

The rate of reaction (1) is binary, $\lambda_1=\kappa_1xy$, but one of the factors, $x$, is the variable for the slow process which has constant rate of increase, hence By Proposition~\ref{prop-binbound}, this model has an exponentially good approximation in a sequence of processes, indexed by increasing values of $M'$, where only the rate of reaction (1) is replaced by $\lambda'_1(z)=\kappa_1y(x\wedge M')$. 
 
The exponential generator for $f_N(x,y)=f(x)+N^{-1}g(x,y)$ is:
\begin{eqnarray*}
H_Nf_N(x,y)&=&\kappa_0(e^{N(f(x+N^{-1})-f(x))}-1)\\
&+&\kappa_1xy(e^{N(f(x-N^{-1})-f(x))+g(x-N^{-1},y-1)-g(x,y)}-1)\\
&+&\kappa_2x(e^{g(x,y+1)-g(x,y)}-1)\\
&+&\kappa_3y(e^{g(x,y-1)-g(x,y)}-1)\\
&+&\kappa_4(e^{g(x,y+1)-g(x,y)}-1)
\end{eqnarray*}
so that its limit $\ol H_0(x,p)$ solves the eigenvalue problem \eqref{eq-EVP} with 
\begin{eqnarray*}
V(y;x,p)&=&\kappa_0(e^{p}-1)+\kappa_1xy(e^{-p}-1)\\
L_1^{x,p}e^{g(x,y)}&=&\kappa_1xye^{-p}(e^{g(x,y-1)}-e^{g})+\kappa_2x(e^{g(x,y+1)}-e^{g})+\kappa_3y(e^{g(x,y-1)}-e^{g})+\kappa_4(e^{g(x,y+1)}-e^{g})
\end{eqnarray*}
Letting $g(x,y)=a(x)y$ we get:
\begin{eqnarray*}
\ol H_0(x,p)&=&\kappa_0(e^{p}-1)+\kappa_1xy(e^{-p-a(x)}-1) +\kappa_2x(e^{a(x)}-1)+\kappa_3y(e^{-a(x)}-1)+\kappa_4(e^{a(x)}-1)\\
&=&y\big(\kappa_1x(e^{-p-a(x)}-1)+\kappa_0(e^{p}-1)+\kappa_3(e^{-a(x)}-1)\big)+\kappa_2x(e^{a(x)}-1)+\kappa_4(e^{a(x)}-1)
\end{eqnarray*}
In order to solve the EVP we set the coefficient of $y$ to $0$ which gives the equation
\[e^{a(x)}=\frac{\kappa_1xe^{-p}+\kappa_3}{\kappa_1x+\kappa_3}>0,\forall x\ge 0, \forall p\]
and hence
\[\ol H_0(x,p)=\kappa_0(e^{p}-1)+\frac{(\kappa_2x+\kappa_4)(\kappa_1xe^{-p}+\kappa_3)}{\kappa_1x+\kappa_3}\]
which is convex and coercive in $p$ and satisfies \eqref{eq:H2ploc} in  $(\ep,\infty)\times \mathbb R$. The effect of the needed truncation will be the replacement of $\kappa_1x$ by $\kappa_1(x\wedge M')$ in the above formula for $\ol H_0$.
\\

To satisfy Condition~\ref{cond-multiplyapun} we note that $|V(y;x,p)|\to\infty$ as $y\to\infty$. For any $c>1$  we let $\varphi_{x,p}=a_{x,p}y$, with $a_{x,p}$ to be chosen. Then calculating as above 
\begin{eqnarray*}
&&e^{-\varphi_{x,p}(y)}L_1^{x,p}e^{\varphi_{x,p}(y)}+c|V(y;x,p)|\\
&&\quad=\kappa_1xye^{-p}(e^{-a_{x,p}}-1)+\kappa_2x(e^{a_{x,p}}-1)+\kappa_3y(e^{-a_{x,p}}-1)+\kappa_4(e^{a_{x,p}}-1)+c\kappa_1xy|e^{-p}-1|\\
&&\quad=y\big((\kappa_1xe^{-p}+\kappa_3)(e^{-a_{x,p}}-1)+c\kappa_1x|e^{-p}-1|\big)+(\kappa_2x+\kappa_4)(e^{a_{x,p}}-1)
\end{eqnarray*}
and chosing $a_{x,p}$ which sets the coefficient of $y$ to 0 
\[e^{a_{x,p}}=\frac{\kappa_1xe^{-p}+\kappa_3}{\kappa_1xe^{-p}+\kappa_3-c\kappa_1x|e^{-p}-1|}\ge 1\]
implies $\varphi_{x,p}=a_{x,p}y\to\infty$  as $y\to\infty$ and
\[e^{-\varphi_{x,p}(y)}L_1^{x,p}e^{\varphi_{x,p}(y)}+c|V(y;x,p)|=(\kappa_2x+\kappa_4)(e^{a_{x,p}}-1),\; \forall y\in E_Y\]
with the right hand side  $d=(\kappa_2x+\kappa_4)(e^{a_{x,p}}-1)\in [0,\infty)$ being independent of $y$, as needed for \eqref{eq-multiplyapun}.

\

\subsection{Viral production (VP)}\label{subsec-VP}
.\\
A somewhat more complicated model is one for production of packaged virus particles:

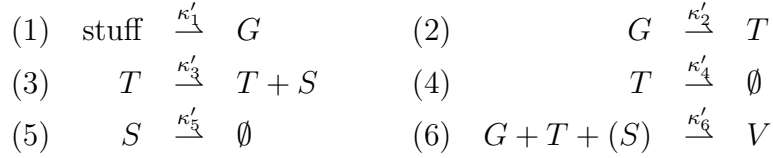
\begin{center}
\begin{tabular}{lrlllrlll}
(1)&$\mbox{\rm stuff}$&$\stackrel {\kappa'_1}{\rightharpoonup}$&$G$&\qquad
(2)&\qquad$G$&$\stackrel {\kappa'_2}{\rightharpoonup}$&$T$\\
(3)&$T$&$\stackrel {\kappa'_3}{\rightharpoonup}$&$T+S$&\qquad
(4)&$T$&$\stackrel {\kappa'_4}{\rightharpoonup}$&$\emptyset$\\
(5)&$S$&$\stackrel {\kappa'_5}{\rightharpoonup}$&$\emptyset$&\qquad
(6)&$G+T+(S)$&$\stackrel {\kappa'_6}{\rightharpoonup}$&$V$\\
\end{tabular}
\end{center}
where $T$ is the viral template, $G$ the viral genome, $S$ the viral structural protein and $V$ the packaged virus. The virus has very few templates from which it manages to co-opt the cell's MRNA to make a relatively large copy number of its genomes, and an order of magnitude larger number of viral structural proteins.  Letting $Z_1,Z_2,Z_3,Z_4$ denote the amounts of species $T,G,S,V$ respectively, the appropriate rescaling gives $Z^N_1=Z_1,Z^N_2=Z_2/N^{2/3},Z^N_3=Z_3/N,Z^N_4=Z_4/N^{2/3}$. The chemical rates also have a range of orders of magnitude, relative to species rescaling they are best expressed by
$\kappa_{1}=\kappa'_{1},\kappa_{2}=\kappa'_{2}N^{2/3}, \kappa_{3}= \kappa'_{3}/N, \kappa_{4}=\kappa'_{4}, \kappa_{5}= \kappa'_{5},\kappa_{6} =\kappa'_{6}N^{5/3}$. The only reaction that is not in standard mass-action form is (6) where the effect of viral proteins are felt only in terms of their order of magnitude, and the usual dependence on the amounts is binary in the amounts of viral templates and genomes: $\lambda_6(z)=\kappa_6z_2z_1$. The model for this system is
\begin{align*}
Z^N_1(t)&=Z^N_1(0)+Y_2(\int_0^t N^{2/3}\kappa_2Z^N_2(s)ds)-Y_4(\int_0^t N^{2/3}\kappa_4Z^N_1(s)ds)\\
&\qquad\qquad\qquad -Y_6(\int_0^t N^{2/3}\kappa_6Z^N_1(s)Z^N_2(s)ds)\\
Z^N_2(t)&=Z^N_2(0)+N^{-2/3}Y_1(\int_0^t N^{2/3}\kappa_1ds)-N^{-2/3}Y_2(\int_
0^t N^{2/3}\kappa_2Z^N_2(s)ds)\\
&\qquad\qquad\qquad -N^{-2/3}Y_6(\int_0^t N^{2/3}\kappa_6Z^N_1(s)Z^N_2(s)ds)\\
Z^N_3(t)&=Z^N_3(0)+N^{-1}Y_3(\int_0^t N^{5/3}\kappa_3Z^N_1(s)ds)-N^{-1}Y_5(\int_
0^t N^{5/3}\kappa_5Z^N_3(s)ds)\\
&\qquad\qquad\qquad  -N^{-1}Y_6(\int_0^t N^{2/3}\kappa_6Z^N_1(s)Z^N_2(s)ds)\\
Z^N_4(t)&=Z^N(0)+N^{-2/3}Y_6(\int_0^t N^{2/3}\kappa_6Z^N_1(s)Z^N_2(s)ds)
\end{align*}
The packaged virus $V$ is the final product but not an ``active'' species in the system, and can be tracked from knowledge on the behaviour of $Z_2^N$.
The slow and fast process are respectively: $X^N=Z^N_2$  and $Y^N=(Z^N_1,Z^N_3)$ and the time-scale separation is now $N^{2/3}$ which will also give the scaling for the large deviation asymptotics as $N^{2/3}$, rather than $N$. 

This is a modified version of the viral production model considered in \cite{BKPR06}, where we showed that the effective dynamics of the fast process $Y^N$ is a piecewise deterministic Markov process with discrete component $Y^N_1$ and continuous component $Y^N_2$. The discrete component is a birth-death Markov chain with birth rate $\kappa_2x$ and death rate $\kappa_4y_1+\kappa_6xy_1$. The continuous component follows the ODE:\;  $\dot{y_2}(t)=\kappa_3y_1-\kappa_5y_2$ which depends on the value $y_1$ of the discrete component $Y^N_1$. This process  has a unique stationary distribution $\mu^x(y_1,y_2)$ which satisfies
\begin{align*}
\int \Big[\kappa_2x\left(g(y_1+1,y_2)-g(y_1,y_2)\right)&+(\kappa_4y_1+\kappa_6xy_1)\left(g(y_1-1,y_2)-g(y_1,y_2)\right)\\&+(\kappa_3y_1-\kappa_5y_2)\partial_{y_2} g(y_1,y_2)\Big]\mu_x(y_1,y_2)=0
\end{align*}
In particular for the discrete component $\mu^x$ has a Poisson($m^x$) distribution with parameter $m^x=\kappa_2x/(\kappa_4+\kappa_6x)$ so $E_{\mu^x}[Y_1]=V_{\mu^x}[Y_1]=m^x$. Moreover, the mean and variance of the continuous component satisfy $E_{\mu^x}[Y_2]=(\kappa_3/\kappa_5)m^x$. 
Using the above results and averaging techniques it can be shown (in the same way as in \cite{BKPR06}) that the effective dynamics of the slow process is given by the ODE 
\[\dot{x}(t)=\kappa_1-\kappa_2 x(t)dt -\kappa_6\frac{\kappa_3}{\kappa_5} \frac{\kappa_2x(t)}{\kappa_4+\kappa_6x(t)}x(t)dt.\]

The rate of reaction (6) is binary, $\lambda_6=\kappa_6xy_1$, but one of the factors is the variable for the slow process which has a constant rate of increase.  By Proposition~\ref{prop-binbound}, this model has an exponentially good approximation in a sequence of processes, indexed by increasing values of $M'$, where only the rate of reaction (6) is replaced by $\lambda'_6(z)=\kappa_6y_2(x\wedge M')$. The effect of the needed truncation will be the replacement of $\kappa_6x$ by $\kappa_6(x\wedge M')$ in the final formula for $\ol H_0$.\\


The generator for the triple is:
\begin{align*}
\vspace{-2mm}
&A_Nf(x,y_1,y_2)=\\
&\quad N\kappa_1\left(f(x+N^{-2/3},y_1,y_2)-f(x,y_1,y_2)\right)
+N^{2/3}\kappa_2x\left(f(x-N^{2/3},y_1+1,y_2)-f(x,y_1,y_2)\right)\\
&+N^{5/3}\kappa_3y_1\left(f(x,y_1,y_2+N^{-1})-f(x,y_1,y_2)\right)
+N^{2/3}\kappa_4y_1\left(f(x,y_1-1,y_2)-f(x,y_1,y_2)\right)\\
&+N^{5/3}\kappa_5y_2\left(f(x,y_1,y_2-N^{-1})-f(x,y_1,y_2)\right)
+N^{2/3}\kappa_6xy_1\left(f(x-N^{2/3},y_1-1,y_2)-f(x,y_1,y_2)\right)
\end{align*}
so the exponential generator, on functions of the form $f_N(x,y)=f(x)+N^{-2/3}g(x,y_1,y_2)$ is:
\begin{align*}
\vspace{-2mm}
&H_Nf_N(x,y_1,y_2)=\\
&\quad \kappa_1\big(e^{N^{2/3}(f(x+N^{-2/3})-f(x))+g(x+N^{-2/3},y_1,y_2)-g(x,y_1,y_2)}-1\big)\\
&+\kappa_2x\big(e^{N^{2/3}(f(x-N^{2/3})-f(x))+g(x-N^{-2/3},y_1+1,y_2)-g(x,y_1,y_2)}-1\big)\\
&+\kappa_3y_1\big(e^{g(x,y_1,y_2+N^{-1})-g(x,y_1,y_2)}-1\big)+\kappa_4y_1\big(e^{g(x,y_1-1,y_2)-g(x,y_1,y_2)}-1\big)+\kappa_5y_2\big(e^{g(x,y_1,y_2-N^{-1})-g(x,y_1,y_2)}-1\big)\\
&+\kappa_6xy_1\big(e^{N^{2/3}(f(x-N^{2/3})-f(x))+g(x,y_1-1,y_2)-g(x,y_1,y_2)}-1\big)
\end{align*}
and its limit is:
\begin{align*}
\vspace{-2mm}
\lim\limits_{N\to\infty}&H_Nf_N(x,y_1,y_2)\\ &=\kappa_1\big(e^{f'(x)}-1\big)+\kappa_2x\big(e^{-f'(x)+g(x,y_1+1,y_2)-g(x,y_1,y_2)}-1\big)+\kappa_3y_1\partial_{y_2}g(x,y_1,y_2)\\
&+\kappa_4y_1\big(e^{g(x,y_1-1,y_2)-g(x,y_1,y_2)}-1\big)-\kappa_5y_2\partial_{y_2}g(x,y_1,y_2)+\kappa_6xy_1\big(e^{-f'(x)+g(x,y_1-1,y_2)-g(x,y_1,y_2)}-1\big)
\end{align*}

We can identify the potential function $V$ as
\begin{eqnarray*}
V(y_1,y_2;x,p)=\kappa_1(e^p-1)+\kappa_2x(e^{-p}-1)+\kappa_6xy_1(e^{-p}-1)\end{eqnarray*}
and perturbed generator  $L_1^{x,p}$ for effective fast process as 
\begin{eqnarray*}
L_1^{x,p}e^{g(x,y_1,y_2)}&=&\kappa_2xe^{-p}(e^{g(x,y_1+1,y_2)}-e^{g(x,y_1,y_2)})+\kappa_4y_1(e^{g(x,y_1-1,y_2)}-e^{g(x,y_1,y_2)})\\
&+&(\kappa_3y_1-\kappa_5y_2)\partial_{y_2} e^{g(x,y_1,y_2)}+\kappa_6xy_1e^{-p}(e^{g(x,y_1-1,y_2)}-e^{g(x,y_1,y_2)})
\end{eqnarray*}
To solve the eigenvalue problem \eqref{eq-EVP} with this $V$ and $L_1^{x,p}$,  we let $g(x,y)=a_1(x)y_1+a_2(x)y_2$ which would imply
\begin{eqnarray*}
\ol H_0(x,p)&=&\kappa_1(e^p-1)+\kappa_2x(e^{-p}-1)+\kappa_6xy_1(e^{-p}-1)\\&+&\kappa_2xe^{-p}(e^{a_1(x)}-1)+\kappa_4y_1(e^{-a_1(x)}-1)+(\kappa_3y_1-\kappa_5y_2)a_2(x)+\kappa_6xy_1e^{-p}(e^{-a_1(x)}-1)
\end{eqnarray*}
In the reactions which have linear rates in $y_1$, reactions (4) and (6) have $\tilde\zeta^Y_{4,1}=-1$ in the discrete variable $y_1$, while reaction (3) changes only the continuous variable $y_2$. Hence we get an equation which is potentially quadratic in $e^{a_1(x)}$, and involves $a_2(x)$ as well. However, Condition~\ref{cond-conversionPDMP} is satisfied since there exist only reaction (5) which has linear rate in $y_2$ and it only changes $y_2$. Thus the equation for $a_2(x)$ can be solved independently of $e^{a_1(x)}$, and the one equation for $e^{a_1(x)}$ potentially quadratic in fact has the quadratic coefficient equal to 0. Setting coefficients of $y_1$ and $y_2$ to zero we get
\[ a_2(x)=0,\; e^{-a_1(x)}=\frac{\kappa_4+\kappa_6x-\kappa_3a_2(x)}{\kappa_4+\kappa_6xe^{-p}}>0\;  \forall x\ge 0, \forall p\]
\begin{eqnarray*}
\ol H_0(x,p)&=\kappa_1(e^p-1)+\kappa_2x(e^{-p}-1)(1+\frac{\kappa_6xe^{-p}}{\kappa_4+\kappa_6x})
\end{eqnarray*}
which is convex and coercive in $p$ and satisfies \eqref{eq:H2ploc} in  $(\ep,\infty)\times \mathbb R$. Recall that one needs to replace $\kappa_6x$ by $\kappa_6(x\wedge M')$ in the above formula for $\ol H_0$.

Since the state space for the fast variables in not compact we need to check the Lyapunov Condition~\ref{cond-multiplyapun}. Note that $|V(y_1,y_2;x,p)|\to\infty$ as $|(y_1,y_2)|\to\infty$. 
Let $y_1^*$ be large enough so that $(\kappa_6xy_1^*+\kappa_2x)|e^{-p}-1|>\kappa_1(e^p-1)$.
For any $c>1$ we let $\varphi_{x,p}=a_{x,p}y_1$, then $\forall  y_1>y_1^*, \forall y_2$ 
\begin{eqnarray*}
&&e^{-\varphi_{x,p}(y)}L_1^{x,p}e^{\varphi_{x,p}(y)}+c|V(y;x,p)|\\
&&\quad=\kappa_2xe^{-p}(e^{a_{x,p}}-1)+\kappa_4y_1(e^{-a_{x,p}}-1)+\kappa_6xy_1e^{-p}(e^{a_{x,p}}-1)\\
&+&c((\kappa_6xy_1^*+\kappa_2x)|e^{-p}-1|-\kappa_1(e^p-1))
\end{eqnarray*}
and chosing $a_{x,p}$ to set the coefficients of $y_1$ to 0
\[e^{-a_{x,p}}=\frac{\kappa_4+\kappa_6x-c\kappa_6x|e^{-p}-1|}{\kappa_4+\kappa_6x}<1\]
implies $\varphi_{x,p}=ay_1\to\infty$ as $|y_1|\to\infty$ and 
\[e^{-\varphi_{x,p}(y)}L_1^{x,p}e^{\varphi_{x,p}(y)}+c|V(y;x,p)|= \kappa_2xe^{-p}(e^{a_{x,p}}-1)+c\kappa_2x)|e^{-p}-1|-c\kappa_1(e^p-1)\]
the right-hand side being independent of $(y_1,y_2)$ as needed for \eqref{eq-multiplyapun}.
%
 \

%

\def\ep{\epsilon}
\def\ol{\overline}
\def\ul{\underline}
\def\l{\left}
\def\r{\right}
\def\p{\partial}
\def\xep{\bar{x}_N}
\def\yep{\bar{y}_N}
\def\sep{\bar{s}_N}
\def\tep{\bar{t}_N}

\setcounter{equation}{0}
\section{Appendix}\label{sec-proof}
We adapt the proof of the large deviation principle for the two time-scale jump-diffusions given in \cite{KP17} to the context of our multiscale Markov chain $Z_N$. Since the steps of the proof  are the same, we give a fairly terse outline and explain carefully the details only where the proof is modified.

\subsection{Proof of the LDP Theorem~\ref{thm-ldp}}\label{subsec-proofthm}

\begin{proof}
Our proof is based on the viscosity solution to the Cauchy problem for each $h\in C_b(E)$ 
\begin{equation}\label{eq-CauchyNN}\p_t u_N=H_Nu_N, \mbox{ in }(0,T]\times E; \quad u_N(0,\cdot)=h(\cdot), \mbox{ in } E\end{equation}
where the non-linear operator is the exponential generator $H_Nf=\frac 1Ne^{-Nf}A_Ne^{Nf}$ for $e^{Nf}\in\mathcal D(A_N)$ given by
\[H_Nf(z)=\frac1N\sum_kN^{\beta_k}\lambda^N_k(z)(e^{N(f(z+N^{-\ul \alpha}\zeta_k^N)-f(z))}-1).\]
The definition of viscosity solutions for these types of non-local partial integro-differential equations (PIDEs)  and their properties were given in \cite{Lenhart}, and various results can be found in \cite{BI08}. In order to establish the convergence of $u_N$ we need to use a family of integro-differential operators and a sequence of viscosity sub- and super-solutions to associated Cauchy problems.

For $\theta\in(0,1),\xi\in C_c^1(E_Y)$, using $\varphi$ satisfying  \eqref{eq-Lyapunov}  of the Lyapunov Condition~\ref{cond-lyapun} we define two sequences (over $N$) of functions: 
\begin{eqnarray*}f_{0,N}(x,y)\!\!\!\!&=&\!\!\!\!f_0(x)+\frac 1Ng_0(y),\\
 g_0(y)\!\!\!\!&=&\!\!\!\!(1-\theta)\xi(y)+\theta\varphi(y), \;
f_0(x)=\phi(x)+\gamma\log(1+x^2);\end{eqnarray*} 
 for some $\gamma>0$ and $\phi\in C_c^1(E_X)$, and  
\begin{eqnarray*}f_{1,N}(x,y)\!\!\!\!&=&\!\!\!\!f_1(x)+\frac 1Ng_1(y)
\\ g_1(y)\!\!\!\!&=&\!\!\!\!(1+\theta)\xi(y)-\theta\varphi(y), \;
f_1(x)=\phi(x)-\gamma\log(1+x^2)\vspace{-5mm}
\end{eqnarray*} 
Then,
\begin{eqnarray*}
&&H_Nf_{0,N}(x,y)\\
&&=\sum_{k:\beta_k=1}\lambda^N_k(z)\big(e^{N(f_0(x+N^{-\underline \alpha}\zeta_k^N)-f_0(x))+(g_0(y+N^{- \underline\alpha}\zeta_k^N)-g_0(y))}-1\big)\\
&&\;+\sum_{k:\beta_k>1}N^{\beta_k-1}\lambda^N_k(z)\big(e^{N(f_0(x+N^{- \underline\alpha}\zeta_k^N)-f_0(x))+(g_0(y+N^{- \underline\alpha}\zeta_k^N)-g_0(y))}
-N N^{-\underline\alpha}\zeta^N_k\cdot\nabla f_0(x)- N^{-\underline\alpha}\zeta^N_k\cdot\nabla g_0(y)-1\big)\\
&&\;+\sum_{k:\beta_k>1}N^{\beta_k}\lambda^N_k(z) N^{-\underline\alpha}\zeta_k\cdot\nabla f_0(x)+\sum_{k:\beta_k>1}N^{\beta_k-1}\lambda^N_k(z)N^{-\underline\alpha}\zeta_k\cdot\nabla g_0(y)\\
&&=\sum_{k:\beta_k=1} \lambda^N_k(z)\big(e^{\tilde\zeta^X_k \cdot \nabla f_0}-1\big)
+\sum_{k:\beta_k=1} \lambda^N_k(z)e^{\tilde\zeta^X_k \cdot \nabla f_0}\big(e^{g_0(y+\tilde\zeta^Y_k)-g_0(y)}-1\big)\\
&&\;+\sum_{k:\beta_k>1}\lambda^N_k(z) \tilde \zeta^X_k \cdot \nabla f_0+\sum_{k:\beta_k>1}\lambda^N_k(z) \tilde \zeta^Y_k\cdot\nabla g_0(y) +N^{-1}\varepsilon_N(x,y)\\
&&=V(y;x,\nabla f_0)+\sum_{k:\beta_k=1} \lambda^N_k(z)e^{\tilde\zeta^X_k \cdot \nabla f_0}\big(e^{g_0(y+\tilde\zeta^Y_k)-g_0(y)}-1\big)+\sum_{k:\beta_k>1}\lambda^N_k(z) \tilde \zeta^Y_k\cdot\nabla g_0(y) +N^{-1}\varepsilon_N(x,y)
\end{eqnarray*}
with $|\varepsilon_N|$ bounded on compact subsets of $E_X\times E_Y$ (note $\sup_x|\nabla f_0(x)|<\infty$). \\
Letting $p=\nabla f_0$ we have (by convexity of $e$)
\begin{eqnarray*}
&&H_Nf_{0,N}(x,y)\\
&&\quad\le V(y;x,p)+\sum_{k:\beta_k=1} \lambda^N_k(z)e^{\tilde\zeta^X_k \cdot p}\big((1-\theta)(e^{\xi(y+\tilde\zeta^Y_k)-\xi(y)}-1)+\theta (e^{\varphi(y+\tilde\zeta^Y_k)-\varphi(y)}-1)\big)\\
&&\quad\qquad\qquad\quad +\sum_{k:\beta_k>1}\lambda^N_k(z)\tilde \zeta^Y_k\cdot\big((1-\theta)\nabla \xi(y)+\theta\nabla \varphi(y)\big) +N^{-1}\varepsilon_N(x,y)\\
&&\quad=  V(y;x,p)+(1-\theta)e^{-\xi}L_1^{x,p}e^{\xi}(y)+\theta e^{-\varphi}L_1^{x,p}e^{\varphi}(y)+N^{-1}\varepsilon_N(x,y)\end{eqnarray*}
so 
\[\limsup_{N\to\infty}H_Nf_{0,N}(x,y)\le V(y;x,p)+(1-\theta)e^{-\xi}L_1^{x,p}e^{\xi}(y)+\theta e^{-\varphi}L_1^{x,p}e^{\varphi}(y).\]
Similarly 
\begin{eqnarray*}
&&H_Nf_{1,N}(x,y)\\
&&\quad= V(y;x,\nabla f_1)+\sum_{k:\beta_k=1} \lambda^N_k(z)e^{\tilde\zeta^X_k \cdot \nabla f_1}\big(e^{g_1(y+\tilde\zeta^Y_k)-g_1(y)}-1\big)+\sum_{k:\beta_k>1}\lambda^N_k(z) \tilde \zeta^Y_k\cdot\nabla g_1(y) +N^{-1}\varepsilon_N(x,y)\\
&&\quad \ge 
V(y;x,p)+(1+\theta)e^{-\xi}L_1^{x,p}e^{\xi}(y)-\theta e^{-\varphi}L_1^{x,p}e^{\varphi}(y)+N^{-1}\varepsilon_N(x,y)
\end{eqnarray*}
where now  $p=\nabla f_1$ (and we used inequalities $e^x-1\ge x\ge 1-e^{-x}$).
So 
\[\liminf_{N\to\infty}H_Nf_{1,N}(x,y)\ge  V(y;x,p)+(1+\theta)e^{-\xi}L_1^{x,p}e^{\xi}(y)-\theta e^{-\varphi}L_1^{x,p}e^{\varphi}(y)\]
This implies that for any $\{x_N,y_N\}$ contained in a compact subset of $E_X\times E_Y$ with $x_N\to x$
\[\limsup_{N\to\infty}H_Nf_{0,N}(x_N,y_N)\le H_0(x,p;\xi,\theta)\] and \[\liminf_{N\to\infty}H_Nf_{1,N}(x_N,y_N)\ge H_1(x,p;\xi,\theta)\] where the Lyapunov Condition~\ref{cond-lyapun} allows us to define  two  families (over \mbox{$\theta\in(0,1),\xi\in C_c^1(E_Y)$}) of operators
\begin{eqnarray*}
&&H_0(x,p;\xi,\theta):=
\sup_{y\in E_Y}\{V(y;x,p) +(1-\theta)e^{-\xi}L_1^{x,p}e^{\xi}(y)+\theta e^{-\varphi}L_1^{x,p}e^{\varphi}(y)\},\\
&&H_1(x,p;\xi,\theta):=\inf_{y\in E_Y}\{V(y;x,p)
+(1+\theta)e^{-\xi}L_1^{x,p}e^{\xi}(y)-\theta e^{-\varphi}L_1^{x,p}e^{\varphi}(y)\}.
\end{eqnarray*}
It is easily seen  (again using Condition~\ref{cond-lyapun} and the fact that $\xi\in C^1_c(E_Y)$) that for $c>0$ \mbox{$\{H_{N,0}f_{0,N}\ge -c\}\cap\{f_{0,N}\le c\}$} and \mbox{$\{H_{N,0}f_{1,N}\le c\}\cap\{f_{1,N}\ge -c\}$} are contained in compact subsets of $E_X\times E_Y$ and  by construction  $\{f_{0,N}\}$ and $\{f_{1,N}\}$ converge uniformly on compact subsets of $E_X\times E_Y$ to $f_0$ and $f_1$ respectively.

The defined sequence of functions thus verifies conditions needed to establish the following result (Condition 3.1 and 3.2 for Lemma 6 in \cite{KP17}): suppose for a sequence of uniformly bounded (over $N$) viscosity solutions $u^h_N$ of the Cauchy problem \eqref{eq-CauchyNN} we construct the upper semicontinuous regularization $\ol u^h$ of the function 
\[u_{\uparrow}^h:=\sup_{y_N}\{\limsup_{N\to\infty}u_N^h(t_N,x_N,y_N):(t_N,x_N)\to (t,x)\}\]
and, similarly,  we construct the lower semicontinuous regularization $\ul u^h$ of the function
\[u_{\downarrow}^h:=\inf_{y_N}\{\liminf_{N\to\infty}u_N^h(t_N,x_N,y_N):(t_N,x_N)\to (t,x)\}\] 
then (see Lemma 6 in \cite{KP17} for details of this conclusion) $\ol u^h$
is a sub-solution of $\p_tu \le H_0(x,\nabla u)$ and $\ul u^h$ is a super-solution of $\p_tu \ge H_1(x,\nabla u)$ with the same initial conditions $u(0,\cdot)=h$, where the two non-linear operators above are defined from the two earlier constructed families of operators by
\[H_0(x,p):=\inf_{\xi,\theta}H_0(x,p;\xi,\theta),\quad H_1(x,p):=\sup_{\xi,\theta}H_1(x,p;\xi,\theta).\]

By definition of $u_{\uparrow}^h$ and $u_{\downarrow}^h$ we immediately have $\ol u^h\ge \ul u^h$,  so if we establish the reverse inequality, we will have proved (Lemma 7 \cite{KP17}) uniform convergence (over compact subsets of $[0,T]\times E_X\times E_Y$) of $u_N$ to $u_0$  the viscosity solution to the Cauchy problem 
\begin{equation}\label{eq-Cauchy00}\p_t u_0=\ol H_0u_0, \mbox{ in }(0,T]\times E; \quad u_0(0,\cdot)=h(\cdot), \mbox{ in } E.\end{equation}

To establish  $\ol u^h\le \ul u^h$ we need to verify that the comparison principle holds between sub-solutions of $\p_t\ol u^h \le H_0(x,\nabla \ol u^h)$ and super-solutions of $\p_t\ul u^h \ge H_1(x,\nabla \ul u^h)$, which will follow if we establish the operator inequality 
\begin{equation}\label{eq-opineq}
H_0(x,p)
\le \ol H_0(x,p) \le  H_1(x,p)
\end{equation}
where $\ol H_0(x,p)$ is defined as the solution to the eigenvalue problem 
\begin{equation}\label{eq-EVPP}
(V(y;x,p)+L_1^{x,p})e^{h(x,y)}=\ol H_0(x,p)e^{h(x,y)}
\end{equation}

To prove the left-hand side we can use the Donsker-Vardhan (see \cite{DV75}) variational representation of the principal eigenvalue $\ol H_0(x,p)$ in \eqref{eq-EVPP} which is given by 
\begin{equation}\label{eq-DVvar}
\ol H_0(x,p)=\sup_{\mu\in \mathcal P(E_Y)}\Big(\int_{E_Y}V(y; x, p)d\mu(y) +\inf_{g\in D^{++}(L_1^{x,p})}\int_{E_Y}\frac{L_1^{x,p}g}{g}(y)d\mu(y) \Big)
\end{equation}
(with $D^{++}(L_1^{x,p})\subset C_b(E_Y)$ denoting functions strictly bounded below by a positive constant). We then repeat the same argument given in Lemma 11.35 of \cite{FK06} that proves
\begin{eqnarray*}&&H_0(x,p)=\inf_{\xi,\theta}\sup_{\mu\in \mathcal P(E_Y)}\int_{E_Y}\Big(V(y;x,p) +(1-\theta)e^{-\xi}L_1^{x,p}e^{\xi}(y)+\theta e^{-\varphi}L_1^{x,p}e^{\varphi}(y)\Big)d\mu(y)\\&&\qquad\qquad \le \sup_{\mu\in \mathcal P(E_Y)}\Big(\int_{E_Y}V(y; x, p)d\mu(y) +\inf_{g\in D^{++}(L_1^{x,p})}\int_{E_Y}\frac{L_1^{x,p}g}{g}(y)d\mu(y) \Big)=\ol H_0(x,p)\end{eqnarray*}
(note that in \cite{FK06} operators we denoted as $H_0(x,p), H_1(x,p)$ are indexed by $H_1(x,p), H_2(x,p)$, respectively). The only difference is that we use the density Condition~\ref{cond-ergod} and the Lyapunov Condition~\ref{cond-lyapun} on $\varphi$  to replace their analogous Condition~ 11.21, and to insure the finiteness of the sum of integrals in the above expression.

To prove the right-hand side we first use Lemma B.10 of \cite{FK06} which proves the inequality 
\begin{eqnarray*}&&H_1(x,p)=\sup_{\xi,\theta}\inf_{y\in E_Y}\{V(y;x,p)+(1+\theta)e^{-\xi}L_1^{x,p}e^{\xi}(y)-\theta e^{-\varphi}L_1^{x,p}e^{\varphi}(y)\}\\&&\qquad\qquad \ge \inf_{\mu\in \mathcal P(E_Y)}\liminf_{t\to\infty}\frac 1t \log E^{\mu}\big[e^{\int_0^t V(Y(s);x,p)ds}\big]\end{eqnarray*}
We next use Lemma 8 of \cite{KP17} which uses the density Condition~\ref{cond-ergod} (and the fact that $\inf_yV(y;x,p)>-\infty$) to insure a uniform (over initial points $Y(0)$) large deviation principle lower bound for the occupation measure of the effective dynamics of the fast process $(Y(s); x,p)$ (given by the generator $L_1^{x,p}$)
to prove that  for any initial $\mu\in\mathcal P(E_Y)$ for $(Y(s); x,p)$ we have the inequality 
\[\liminf_{t\to\infty}\frac 1t \log E^{\mu}\big[e^{\int_0^t V(Y(s);x,p)ds}\big]\ge \ol H_0(x,p)\]
with $\ol H_0(x,p)$ given by the variational form \eqref{eq-DVvar} above. 

Having proved the Operator Inequality \eqref{eq-opineq} the constructed sub- and super-solutions sandwich the sequence of viscosity solutions $u_N^h$, establishing their convergence to the viscosity solution $u_0^h$, which is by assumption on the limiting Cauchy problem unique. The next step is to prove   exponential tightness of the sequence $\{X_N(t)\}$. Exponential tightness of paths of $\{X_N\}$  follows from the convergence of the exponential generators $H_N$ by a standard argument (see Lemma 2 in \cite{KP17}, or equivalently see Corollary 4.17 and a simple calculation from Lemma 4.22 in \cite{FK06}). 
Finally, by Bryc formula (Theorem~\ref{thm-varbryc}) we have established a large deviation
principle for $\{X_N(t)\}$ with speed $1/N$ and good rate function $I$ given by $u^h_0(t)$ in terms of a variational principle \eqref{eq-rate}. 
\end{proof}

\

We now use the solution to the eigenvalue problem \eqref{eq-EVP} to simplify the proof of the above theorem by using the associated eigenfunction in forming the family of operators in the definitions of $H_0$ and $H_1$. We also replace the Lyapunov Condition~\ref{cond-lyapun}
by the less stringent Condition~\ref{cond-multiplyapun} (sufficient for proving multiplicative ergodicity of Markov processes, \cite{KM05} of single scale processes). Since  for many chemical reaction models  one can explicitly solve for $\ol H_0(x,p)$ as well as verify \eqref{eq-multiplyapun}, this result is used in all the examples in this paper. 

\subsection{Proof of the LDP Corollary~\ref{cor-ldp}}\label{subsec-proofcor}

\begin{proof}
For each $x\in E_X, p\in K\subset \mathbb R$ let $e^{\xi_{x,p}}\in \mathcal D(\bar H_0)$ denote the eigenfunction associated with the eigenvalue problem \eqref{eq-EVP}, that is, \[V(y;x,p) + e^{-\xi_{x,p}(y)}L^{x,p}_1e^{\xi_{x,p}(y)}=\bar H_0(x,p), \;\forall y E_Y.\]  
For $\theta\in(0,1)$ we let  $\varphi_{x,p}$ be the function satisfying \eqref{eq-multiplyapun}  of the Lyapunov Condition~\ref{cond-lyapun}, and define the following two sequences of functions: 
\begin{eqnarray*}f_{0,N}(x,y)\!\!\!\!&=&\!\!\!\!f_0(x)+\frac 1N((1-\theta)\xi_{x,p}(y)+\theta\varphi_{x,p}(y)),\\
 f_0(x)\!\!\!\!&=&\!\!\!\!\phi(x)+\gamma\log(1+x^2)\end{eqnarray*} 
 for some $\gamma>0$ and $\phi\in C_c^1(E_X)$, as well as  
\begin{eqnarray*}f_{1,N}(x,y)\!\!\!\!&=&\!\!\!\!f_1(x)+\frac 1N((1+\theta)\xi_{x,p}(y)-\theta\varphi_{x,p}(y))\\ 
f_1(x)\!\!\!\!&=&\!\!\!\!\phi(x)-\gamma\log(1+x^2), 
\vspace{-5mm}
\end{eqnarray*} 
Then, for $p=\nabla f_0$, we have (using eigenfunction property of $\xi_{x,p}$ and \eqref{eq-multiplyapun} for $\varphi_{x,p}$)
\begin{eqnarray*}
H_Nf_{0,N}(x,y)\!\!\!\!\!\!\!\!
&&=V(y;x,p)+\sum_{k:\beta_k=1} \lambda^N_k(z)e^{\tilde\zeta^X_k \cdot p}\big(e^{((1-\theta)\xi_{x,p}+\theta\varphi_{x,p})(y+\tilde\zeta^Y_k)-((1-\theta)\xi_{x,p}+\theta\varphi_{x,p})(y)}-1\big)\\
&&\qquad\qquad\quad\;+\sum_{k:\beta_k>1}\lambda^N_k(z) \tilde \zeta^Y_k\cdot\nabla ((1-\theta)\xi_{x,p}+\theta\varphi_{x,p})(y)+N^{-1}\varepsilon(x,y)\\
&&\le  V(y;x,p)+(1-\theta)e^{-\xi_{x,p}}L_1^{x,p}e^{\xi_{x,p}}(y)+\theta e^{-\varphi_{x,p}}L_1^{x,p}e^{\varphi_{x,p}}(y)+N^{-1}\varepsilon(x,y)\\
&&=(1-\theta)\bar H_0(x,p)+\theta V(y;x,p)+\theta e^{-\varphi_{x,p}}L_1^{x,p}e^{\varphi_{x,p}}(y)+N^{-1}\varepsilon(x,y)\\
&& \le (1-\theta)\bar H_0(x,p)-\theta \big((c-1)V(y;x,p)+d\big)+N^{-1}\varepsilon(x,y)
\end{eqnarray*}
Hence, for $c'>0$ we have that \mbox{$\{H_{N,0}f_{0,N}\ge -c'\}\cap\{f_{0,N}\le c'\}$} is contained in compact subset of $E_X\times E_Y$ (since $V(y;x,p)$ has compact level sets and $\bar H_0$ is finite), and for any $\{x_N,y_N\}$ in a compact subset of $E_X\times E_Y$ with $x_N\to x$
\[ \limsup_{N\to\infty} H_Nf_{0,N}(x_N,y_N)\le \sup_{y\in E_Y}\big\{(1-\theta)\bar H_0(x,p)-\theta \big((c-1)V(y;x,p)+d\big)\big\}=:H_0(x,p;\theta)\]
(which by \eqref{eq-multiplyapun} is well defined). 
Similarly, for $p=\nabla f_1$,
\begin{eqnarray*}
H_Nf_{0,N}(x,y)\!\!\!\!\!\!\!\!
&&\ge  V(y;x,p)+(1+\theta)e^{-\xi}L_1^{x,p}e^{\xi}(y)-\theta e^{-\varphi}L_1^{x,p}e^{\varphi}(y)+N^{-1}\varepsilon(x,y)\\
&&=  (1+\theta)\bar H_0(x,p)-\theta V(y;x,p)-\theta e^{-\varphi}L_1^{x,p}e^{\varphi}(y)+N^{-1}\varepsilon(x,y)\\
&&\ge (1+\theta) \bar H_0(x,p)+\theta \big((c-1)V(y;x,p)+d\big)+N^{-1}\varepsilon(x,y),
\end{eqnarray*}
for $c'>0$  \mbox{$\{H_{N,1}f_{1,N}\le c'\}\cap\{f_{1,N}\ge -c'\}$} is contained in compact subsets of $E_X\times E_Y$, and for any $\{x_N,y_N\}$ in a compact subset of $E_X\times E_Y$ with $x_N\to x$
\[ \limsup_{N\to\infty} H_Nf_{1,N}(x_N,y_N)\ge H_1(x,p;\theta):=\inf_{y\in E_Y}\big\{(1+\theta)\bar H_0(x,p)+\theta \big((c-1)V(y;x,p)+d\big)\big\}\]
As all the conditions of Lemma 6 in \cite{KP17} are met, it insures that the two functions $\overline u^h$ and $\underline u^h$ constructed as in the proof of Theorem~\ref{thm-ldp} from the sequence of viscosity solutions of $\partial_t u_N=H_Nu$  are, respectively, a sub-solution of $\partial_t u=H_0 (x,\nabla u)$ and  a super-solution of $\partial_t u=H_1(x,\nabla u)$ with the same initial condition. By construction $\overline u^h\ge \underline u^h$ and the reverse is immediate once we notice that in this case the operators $H_0(x,p)$ and $H_1(x,p)$ given by 
\[H_0(x,p):=\inf_{\theta\in(0,1)}\sup_{y\in E_Y}\big\{(1-\theta)\bar H_0(x,p)-\theta \big((c-1)V(y;x,p)+d\big)\big\}\]
and \[ H_1(x,p):=\sup_{\theta\in(0,1)}\inf_{y\in E_Y}\big\{(1+\theta)\bar H_0(x,p)+\theta \big((c-1)V(y;x,p)+d\big)\big\}\]
 coincide and are equal to $\bar H_0$ (see also \cite{FK06} Appendix B, Lemma 11.4).
\end{proof}

\subsection{Proof of the comparison principle Lemma~\ref{lem-H0cond}}\label{subsec-prooflem}

\begin{lemma} 
Suppose $u_1$ and $u_2$ are,respectively, a bounded upper semicontinuous (USC) viscosity sub-solution and a bounded lower semicontinuous (LSC) viscosity  super-solution of  \eqref{eq-Cauchy0}
for some $T>0$ and $E_X\subset \R^d$. Either of the following conditions are sufficient for the weak comparison principle for \eqref{eq-Cauchy0} to hold:
\begin{itemize}
\item[(a)]  $\ol H_0$ is such that for all $\lambda\ge 1, R>0$ and for all $|p|,|q|\le 1$, $|x|,|y|\le R$ and for some continuous non-decreasing functions $\omega_R, \,\tilde \omega_1:\R_+\mapsto\R_+$ with $\omega_R(0)=\tilde \omega_1(0)=0$
\begin{equation}\label{eq:H2}\begin{aligned}
\hspace{-5mm}\ol H_0&(x, \lambda(x-y)+p)-\ol H_0(y,\lambda(x-y)+q)\le \omega_R(|x-y|+\lambda |x-y|^2)+\tilde\omega_1(|p-q|)\end{aligned}
\end{equation}
\item[(b)]  $\ol H_0$ is such that for all $\lambda\ge 1, R,\ell>0$ and for all $|p|,|q|\le 1$, $|x|,|y|\le R$ \underline{with $\lambda|x-y|<\ell$}, and for continuous non-decreasing functions $\omega_{R,\ell}, \,\tilde \omega_1:\R_+\mapsto\R_+$ with \mbox{$\omega_{R,\ell}(0)=\tilde \omega_1(0)=0$} the inequality \eqref{eq:H2} holds; 
and 
$\ol H_0$ satisfies a coercivity condition in $p$
\begin{equation}\label{eq:coercive}
\frac{\ol H_0(x,p)}{|p|}\to c_x>1\mbox{ as }|p|\to\infty \mbox{ uniformly with respect to $x$},
\end{equation} 
\end{itemize}
\end{lemma}

\begin{proof} \indent {(a)} We first prove the comparison principle under the assumptions on $\ol H_0$ given in (a).
Let $A=\sup_{\{0\}\times E_X}(u_1-u_2)^+$, and suppose  for some $\tilde x$ and $\tilde t$, $u_1(\tilde t,\tilde x)-u_2(\tilde t,\tilde x)=A+\delta$ ($\ast$) with $\delta>0$. 
For $\beta, m, \eta>0, N<\infty$  define \[\psi_N(t,x,s,y)=u_1(t,x)-u_2(s,y)-\frac{N}{2}\l[|x-y|^2+|t-s|^2 \r]-\beta\l(g(x)^m+g(y)^m\r) -\eta(t+s)\]
where $g(x)=\sqrt{1+|x|^2}$.
Now choose $\beta,\eta$ such that $2\beta g(\tilde x)^m+2\eta g(\tilde t)\le \frac\delta 2$ holds for all $m\le 1$. 
Then ($\ast$) implies
\begin{equation}\label{eq:contradiction}
\sup_{[0,T]\times E_X\times [0,T]\times E_X} \psi_N(t,x,s,y)\ge \psi_N(\tilde t, \tilde x,\tilde t,\tilde x)=A+\delta-2\beta g(\tilde x)^m-2\eta \tilde t\ge A+\frac \delta 2
\end{equation}
and we show this gives a contradiction.

Since $\psi_N$ is USC and tends to $ -\infty$ as $|x|+|y|\to\infty$ its maximum in $[0,T]\times E_X\times[0,T]\times E_X$ is achieved at a point $(\tep,\xep,\sep,\yep)$. By \eqref{eq:contradiction} $\psi_N(\tep,\xep,\sep,\yep)\ge A+\frac\delta 2$ and $u_1, u_2$ bounded imply 
$\beta(g(\xep)^m+g(\yep)^m)\le \sup (u_1- u_2)-A-\frac\delta 2=:c_1\ge 2\beta$ for all $N$, all $m\le 1$. Hence, $g(\xep),g(\yep)\le (\frac{c_1}{\beta})^{1/m}$ and for $R:=(\frac{c_1}{\beta})^{1/m}>0$ we have 
\begin{equation}\label{eq:compactmax}
|\xep|,|\yep|\leq R,\; \forall N<\infty.
\end{equation}

Using $\psi_N(\tep,\xep,\tep,\xep)+\psi_N(\sep,\yep,\sep,\yep)\le  2\psi_N(\tep,\xep,\sep,\yep)$ we have \begin{equation*}\frac{N}{2}\l[|\xep-\yep|^2+|\tep-\sep|^2 \r]\le \frac{1}{2}\l[u_1(\tep, \xep)-u_1(\sep,\yep)+u_2(\tep, \xep)-u_2(\sep,\yep)\r]\le c_2\end{equation*}
where $c_2:=\sup |u_1|+\sup |u_2|$. Therefore, $|\xep-\yep|+ |\tep-\sep|\le \sqrt{\frac{2}{N} c_2}$, and $|\xep-\yep|, |\tep-\sep|\to 0$ as $N\to \infty$. We can also show $\frac{N}{2}\l[|\xep-\yep|^2+|\tep-\sep|^2 \r]\to 0$. Let \[S=\max_{|x|\le R,t\le T}\l[(u_1-u_2)(x)- 2\beta g(x)^m-2\eta t\r],\] then for all $N$ due to \eqref{eq:compactmax} we have
\begin{eqnarray*}\begin{aligned}
S&=\max_{|x|\le R,t\le T}\psi_N(t,x,t,x)\le \max_{|x|,|y|\le R,t,s\le T}\psi_N(t,x,s,y)=\psi_N(\tep,\xep,\sep,\yep)\\
&\le u_1(\tep,\xep)-u_2(\sep,\yep)-\beta(g(\xep)^m+g(\yep)^m)-\eta(\tep+\sep):=S_N.
\end{aligned}\end{eqnarray*} 
If we show $\lim_{N\to \infty}S_N\le S$ this will imply that $\lim_{N\to \infty}\frac{N}{2}\l[|\xep-\yep|^2+|\tep-\sep|^2 \r]=0$.  Suppose there exists $N_k\to \infty$ such that $\lim_{N_k\to \infty}S_{N_k}>S$. Since $|\xep|,|\yep|\le R$ and $ \tep,\sep\le T$ we can assume  that $\lim_{N_k\to \infty}(\bar x_{N_k},\bar y_{N_k},\bar t_{N_k},\bar s_{N_k})= (\bar x,\bar y,\bar t, \bar s)$ with $|\bar x|,|\bar y|\le R$, and $ \bar t,\bar s\le T$; and since $|\xep-\yep|+|\tep-\sep|\to 0$ we have $\bar x=\bar y$ and $\bar t=\bar s$. Since $u_1-u_2$ is USC we have 
\begin{eqnarray*}\begin{aligned}
\lim_{N_k\to \infty}&\psi_N(\bar t_{N_k},\bar x_{N_k},\bar s_{N_k},\bar y_{N_k})\le \lim_{N_k\to \infty}S_{N_k}\\&=\lim_{N_k\to \infty}u_1(\bar t_{N_k},\bar x_{N_k})-u_2(\bar s_{N_k},\bar y_{N_k})-\beta ( g(\bar x_{N_k})^m+g(\bar y_{N_k})^m)-\eta (\bar t_{N_k}+ \bar s_{N_k}) \\&\le u_1(\bar t,\bar x)-u_2(\bar t,\bar x)-2\beta g(\bar x)^m-2\eta \bar t \le S
\end{aligned}\end{eqnarray*} 
by definition of $S$, which contradicts the assumption that $\lim_{N_k\to \infty}S_{N_k}>S$.
 
 We next show that for some $\bar N>0$ we have $\tep,\sep>0$  for all $N>\bar N$. Suppose, on the contrary, there exists $N_k\to 0$ such that either $t_{N_k}=0$ or $s_{N_k}=0$ for all $k$. By the same arguments as before  we can assume that $\lim_{N_k\to \infty}(\bar x_{N_k},\bar y_{N_k}, \bar t_{N_k}, \bar s_{N_k})= (\bar x, \bar y, \bar t, \bar s)$ where $\bar x=\bar y$ and $\bar t=\bar s=0$. Then, if $\bar s_{N_k}=0$,
 \[\psi_N(\bar t_{N_k},\bar x_{N_k},\bar s_{N_k},\bar y_{N_k})\le u_1(\bar t_{N_k},\bar x_{N_k})- u_2(0,\bar y_{N_k}),\]
 while if $\bar t_{N_k}=0$, 
  \[\psi_N(\bar t_{N_k},\bar x_{N_k},\bar s_{N_k},\bar y_{N_k})\le u_1(0,\bar x_{N_k})  -u_2(\bar s_{N_k},\bar y_{N_k}).\]
Using the fact that  $u_1-u_2$ is USC and the definition of A we have 
 \[\lim_{N_k\to \infty}\psi_N(\bar t_{N_k},\bar x_{N_k},\bar s_{N_k},\bar y_{N_k})\le u_1(0,\bar x)-u_2(0,\bar x)\le A\]
which contradicts ($\ast$) according to which 
$\psi_N(\tep,\xep,\sep,\yep)  \ge A+\frac \delta 2$ for all $N$ with $\delta>0$.

We now define two test functions $\varphi_1,\varphi_2\in\mathcal C^1([0,T]\times E_X)$ 
 \begin{eqnarray*}\begin{aligned}
\varphi_1(t,x)&:=u_2(\sep,\yep)+\frac{N}{2}\l[|x-\yep|^2+|t-\sep|^2\r]+\beta\l(g(x)^m+g(\yep)^m\r) +\eta(t+\sep)\\
\varphi_2(s,y)&:= u_1(\tep,\xep)-\frac{N}{2}\l[|\xep-y|^2-|\tep-s|^2\r]-\beta\l(g(\xep)^m+g(y)^m\r) -\eta(\tep+s).
\end{aligned}\end{eqnarray*} 
so that $(\tep,\xep)$ is a point of maximum of $u_1(t,x)-\varphi_1(t,x)$ and $(\sep,\yep)$ is a point of minimum of $u_2(s,y)-\varphi_2(s,y)$. At the extremum their derivatives in time are
\[\p_t\varphi_1(\tep,\xep)=N(\tep-\sep)+\eta, \; \p_t\varphi_2(s,y)=N(\tep-\sep)-\eta,\] 
and in space are
\[D_x\varphi_1(\tep,\xep)=N(\xep-\yep)+\gamma \xep,\; D_y\varphi_2(\sep,\yep)=N(\xep-\yep)-\tau \yep,\]  with $\gamma=m\beta g(\xep)^{m-2}, \tau=m\beta g(\yep)^{m-2}$. 
Since $\forall N>\bar N$ we have $\tep,\sep\in(0,T]$, and  $u_1$ and $u_2$ are sub- and super- solutions of $\p_t u_0-\ol H_0(x,D_xu_0)=0$ on $(0,T]\times E_X$, we have
\begin{eqnarray*}\begin{aligned}\label{app:subsuper}
\p_t\varphi_1(\tep,\xep) -\ol{H}_0\l(\xep, D_x\varphi_1(\tep,\xep)\r)\le 0,\\
\p_t\varphi_2(\sep,\yep)-\ol{H}_0\l(\yep,D_y\varphi_2(\sep,\yep)\r)\geq 0.
\end{aligned}\end{eqnarray*}
so that
\begin{eqnarray}\begin{aligned}\label{eq:barH0diff}
2\eta&=\p_t\varphi_1(\tep,\xep)-\p_t\varphi_2(\sep,\yep)\\&\leq \ol{H}_0\l(\xep, N(\xep-\yep)+\gamma \xep\r)- \ol{H}_0\l(\yep, N(\xep-\yep)-\tau \yep\r).
\end{aligned}\end{eqnarray}
%
Using \eqref{eq:H2} letting $\lambda=N$, $p=\gamma\xep$, $q= -\tau \yep$ 
we get that the right-hand side is bounded above by 
\begin{eqnarray*}\begin{aligned}
2\eta &\le \omega_R\Big(|\xep-\yep|+N|\xep-\yep|^2\Big)+\tilde\omega_1\Big(m\beta\l[\xep g(\xep)^{m-2}+ \yep g(\yep)^{m-2}\r]\Big)\\&\le \omega_R\Big(|\xep-\yep|+N|\xep-\yep|^2\Big)+\tilde \omega_1(mc_1).
\end{aligned}\end{eqnarray*}
Since $\beta(g(\xep)^m+g(\yep)^m)\le c_1$ now chosing $m\le \min(1, \frac1{c_1})$ so that $\tilde \omega_1(mc_1)<\eta$, and  taking $N$ large enough so that $\omega_R(|\xep-\yep|+\frac{N}{2}|\xep-\yep|^2)<\eta$, leads to a contradiction in the above inequality. This establishes the comparison principle under the conditions given in (a).\\

{(b)} We now extend the comparison principle under the more relaxed assumptions on $\ol H_0$ given in (b). For many examples the condition \eqref{eq:H2}  is too restrictive and needs to be extended to a more local condition in $p$. This can be done if \eqref{eq-Cauchy0} has a bounded Lipschitz continuous viscosity sub- or super-solution, as we show next.

Assume that $v$ is a bounded Lipschitz continuous viscosity solution of \eqref{eq-Cauchy0}. 
Let $A=\sup_{\{0\}\times E_X}(u_1-v)^+$, and suppose that for some $\tilde x, \tilde t$, $u_1(\tilde t,\tilde x)-v(\tilde t,\tilde x)=A+\delta$ ($\ast$) with $\delta>0$. Since $v$ is a super-solution, the same steps as in the proof of (a) with $v$ replacing $u_2$ hold until and including the inequality \eqref{eq:barH0diff}. 
Moreover,  $\psi_N(\tep,\xep,\sep,\xep)+\psi_N(\tep,\xep,\sep,\xep)\le 2 \psi_N(\tep,\xep,\sep,\yep)$ implies 
\begin{eqnarray*}\begin{aligned}
\frac{N}{2}&|\xep-\yep|^2+\frac{N}{2}|\sep-\tep|^2\\&\le \beta(g(\xep)^m-g(\yep)^m)+(v(\sep,\yep)-v(\sep,\xep)+v(\sep,\yep)-v(\tep,\yep))
\end{aligned}\end{eqnarray*}
hence we have $N|(\tep,\xep)-(\sep,\yep)|\le \beta+\ell_v=:\ell$  
where $\ell_v$ is the Lipschitz constant of $v$. 
Using \eqref{eq:H2} in inequality \eqref{eq:barH0diff} with the same choice of $\lambda=N,p=\gamma\xep,q=-\tau\yep$ as earlier gives
\begin{eqnarray*}\begin{aligned}
2\eta &\le \omega_{R,\ell}\Big(|\xep-\yep|+N|\xep-\yep|^2\Big)+\tilde\omega_1\Big(m\beta\l[\xep g(\xep)^{m-2}+ \yep g(\yep)^{m-2}\r]\Big)\\&\le \omega_{R,\ell}\Big(|\xep-\yep|+N|\xep-\yep|^2\Big)+\tilde \omega_1(mc_1).
\end{aligned}\end{eqnarray*}
which also leads to a contradiction in the above inequality, thus establishing $u_1\le v$ on $[0,T]\times E_X$.
Similarly, letting $A^*=\sup_{\{0\}\times E_X}(v-u_2)^+$ and replacing $u_1$ by $v$ in the steps of the proof of (a) will establish \eqref{eq:barH0diff} for the new set of extrema $(\tep^*,\xep^*,\sep^*,\yep^*)$, as well as $N|\xep^*-\yep^*|\le \beta+\ell_v=:\ell$. Using 
\eqref{eq:H2} for $N|\xep^*-\yep^*|\le\ell$ will lead to a contradiction in \eqref{eq:barH0diff} and establish $v\le u_2$ on $[0,T]\times  E_X$.

It is easy to modify the above steps to show the comparison principle if we have a bounded LSC viscosity sub-solution  $u_1$ and a USC viscosity super-solution $u_2$ of \eqref{eq-Cauchy0} and $\ol  H_0$ satisfying 
\eqref{eq:H2} for $N|\xep^*-\yep^*|\le\ell$, such that either $u_1$ or $u_2$ is Lipschitz continuous. 
To insure Lipschitz continuity of the USC viscosity super-solution $u_2$ we use Lemma~\ref{lem:coerciveLip}, which then establishes the comparison principle under the conditions given in  (b).
\end{proof}

\begin{lemma}\label{lem:coerciveLip} Suppose $\ol H_0$ satisfies the coercivity condition
\begin{equation}\label{eq:coercive}
\frac{\ol H_0(x,p)}{|p|}\to c_x>1\mbox{ as }|p|\to\infty
\end{equation} 
uniformly with respect to $x$ $(\exists C$ such that $\forall x, \forall |p|>C: \ol H_0(x,p)>c_x|p|)$.  If $u_2$ is a bounded lower semicontinuous viscosity super-solution to \eqref{eq-Cauchy0}, then $u_2$ is Lipschitz continuous.
\end{lemma}
\begin{proof} For fixed $(t,x)\in [0,T]\times E_X$ let 
\[\varphi(s,y)=-C\l(|y-x|+|s-t|\r)-u_2(s,y)\] where $C>0$. Since $u_2$ is bounded and LSC there exists $(\bar s,\bar y)\in [0,T]\times E_X$ such that $\varphi(\bar s,\bar y)=\sup_{[0,T]\times  E_X}\varphi(t,y)$. We show that $(\bar s,\bar y)=(t,x)$ for $C$ large enough. Suppose otherwise, then $\varphi+u_2$ is differentiable at all $(s,y)= (\bar s,\bar y)$ such that $\bar s\neq t$ and  $\bar y\neq x$ and since $u_2$ is a super-solution of \eqref{eq-Cauchy0} we would have
\[C\frac{t-\bar s}{|t-\bar s|}-\ol H_0(\bar y,C\frac{x-\bar y}{|x-\bar y|})\ge 0\] 
and for sufficiently large $C$ independent of $x$ this is in contradiction with \eqref{eq:coercive}.
Thus, for such $C$ we have 
\[u_2(t,x)=u_2(\bar s,\bar y)+C\l(|\bar y-x|+|\bar s-t|\r)\le u_2(s,y)+C\l(|y-x|+|s-t|\r)\]
and interchanging the roles of $(t,x)$ and $(s,y)$ shows $u_1$ is Lipschitz continuous.

\end{proof}

\begin{remark} When $\ol H_0$ is strictly convex and in fact super-linear in $p$ one can use calculus of variations and results in deterministic optimal control to show that the value function is the unique Lipschitz continuous viscosity solutions of  the Cauchy problem \eqref{eq-Cauchy0}. The {Fenchel-Legendre transform} provides a formula
\begin{equation}\label{eq:duality}
L(x,v)=\max_{p\in\R^d}\l(-v\cdot p -\ol H_0(x,p)\r)
\end{equation}
for the running cost in a deterministic optimal control problem on a fixed time interval with an initial cost function. Suppose $\ol H_0$ satisfies the following set of conditions:
\begin{eqnarray}\label{eq:Lipcond}\begin{aligned}
(a)\; & \p_{pp}\ol H_0(x,p)>0 \mbox{ (positive definite)},\\
(b)\; & \lim_{|p|\to \infty} \frac{\ol H_0(x,p)}{|p|}=+\infty,\\
(c)\; & p\cdot \p_p\ol H_0  - \ol H_0\ge |\p_p \ol H_0|\gamma(\p_p \ol H_0), \mbox{ where } \gamma(p)\mathop{\to}\limits_{|p|\to\infty} \infty\\
(d)\; & \ol H_0(x,0)\le c_1,\;\; \ol H_0(x,p)\ge -c_1,\\
(e)\; & |\p_x \ol H_0| \le c_2(p\cdot \p_p\ol H_0-\ol H_0)+c_3,\;\;
|\p_p \ol H_0|\le R \Rightarrow |p|\le C(R).
\end{aligned}\end{eqnarray}
The following is Theorem 10.3 from \cite{FS06} Chapter II. 
The conditions in \eqref{eq:Lipcond} can be verified to hold for all the models in our Examples using explicit formulae for $\ol H_0$
\begin{lemma} Let $u_0$ be the value function of the calculus of variation problem on a fixed time interval $[0,T]$ with no constraints on the terminal state,  whose cost function is obtained from $\ol H_0$ by the Fenchel-Legendre duality formula \eqref{eq:duality} and whose initial cost is $f$. If $\ol H_0$ satisfies \eqref{eq:Lipcond}, then $u_0$ is the unique bounded Lipschitz continuous solution to \eqref{eq-Cauchy0}.
\end{lemma}

\end{remark}

\noindent 
{\bf Acknowledgments}. The author would like to thank an anonymous referee for pointing out the incorrect proof of Lemma 3.7 in an earlier version, as well as Rohini Kumar and Tom Kurtz for useful discussions about the material.

\bibliography{chemldp}

\begin{thebibliography}{14}
\providecommand{\natexlab}[1]{#1}
\providecommand{\url}[1]{\texttt{#1}}
\expandafter\ifx\csname urlstyle\endcsname\relax
  \providecommand{\doi}[1]{doi: #1}\else
  \providecommand{\doi}{doi: \begingroup \urlstyle{rm}\Url}\fi

\bibitem[Bakhtin and Hurth(2012)]{BH12}
Yuri Bakhtin and Tobias Hurth.
\newblock Invariant densities for dynamical systems with random switching.
\newblock \emph{Nonlinearity}, 25\penalty0 (10):\penalty0 2937, 2012.

\bibitem[Ball et~al.(2006)Ball, Kurtz, Popovic, and Rempala]{BKPR06}
Karen Ball, Thomas~G. Kurtz, Lea Popovic, and Greg Rempala.
\newblock Asymptotic analysis of multiscale approximations to reaction
  networks.
\newblock \emph{Ann. Appl. Probab.}, 16\penalty0 (4):\penalty0 1925--1961,
  2006.
\newblock ISSN 1050-5164.

\bibitem[Bardi and Dolcetta(1997)]{BD}
Martion Bardi and Italo Capuzzo-Dolcetta.
\newblock \emph{Optimal Control and Viscosity Solutions of HamiltonÐJacobi-Bellman Equations}.
\newblock  Systems \& Control: Foundations \& Applications, Birkh¬auser, Boston, 1997.

\bibitem[Barles and Imbert(2008)]{BI08}
Guy Barles and Cyril Imbert. 
\newblock Second-order elliptic integro-differential equations: viscosity solutions' theory revisited.
\newblock \emph{Annales de l'Institut Henri Poincare (C)  Non Linear Analysis.}
\newblock Vol. 25. No. 3. Elsevier Masson, 2008.

\bibitem[Cloez et~al.(2015)Cloez, Hairer, et~al.]{CH15}
Bertrand Cloez, Martin Hairer, et~al.
\newblock Exponential ergodicity for markov processes with random switching.
\newblock \emph{Bernoulli}, 21\penalty0 (1):\penalty0 505--536, 2015.

\bibitem[Crandall et~al.(1992)Crandall, Ishii, and Lions]{CIL92}
Michael~G Crandall, Hitoshi Ishii, and Pierre-Louis Lions.
\newblock UserÕs guide to viscosity solutions of second order partial
  differential equations.
\newblock \emph{Bulletin of the American Mathematical Society}, 27\penalty0
  (1):\penalty0 1--67, 1992.

\bibitem[Darden(1979)]{Dard79}
Thomas Darden.
\newblock A pseudo-steady state approximation for stochastic chemical kinetics.
\newblock \emph{Rocky Mountain J. Math.}, 9\penalty0 (1):\penalty0 51--71,
  1979.
\newblock ISSN 0035-7596.
\newblock Conference on Deterministic Differential Equations and Stochastic
  Processes Models for Biological Systems (San Cristobal, N.M., 1977).

\bibitem[Davis(1993)]{Dav93}
M.~H.~A. Davis.
\newblock \emph{Markov models and optimization}, volume~49 of \emph{Monographs
  on Statistics and Applied Probability}.
\newblock Chapman \& Hall, London, 1993.
\newblock ISBN 0-412-31410-X.

\bibitem[Dembo and Zeitouni(1998)]{DZ}
Amir Dembo and Ofer Zeitouni.
\newblock Large deviations techniques and applications, 1998.

\bibitem[Donsker and Varadhan(1975)]{DV75}
Monroe~D.~Donsker and S.R.~ Srinivasa Varadhan. 
\newblock On a variational formula for the principal eigenvalue for operators with maximum principle. \newblock \emph{Proceedings of the National Academy of Sciences} 72.3: 780-783  (1975).

\bibitem[Donsker and Varadhan(1983)]{DV83}
Monroe~D.~Donsker and S.R.~ Srinivasa Varadhan. 
\newblock Asymptotic evaluation of certain Markov process expectations for large time. IV
\newblock \emph{Communications on Pure and Applied Mathematics} 36.2:183-212 (1983).


\bibitem[Feng(1999)]{F99}
Jin Feng. 
\newblock Martingale problems for large deviations of Markov processes. 
\newblock \emph{Stochastic Processes and their Applications}, 81.2:165Ð216  (1999).

\bibitem[Feng and Kurtz(2006)]{FK06}
Jin Feng and Thomas~G. Kurtz.
\newblock \emph{Large deviations for stochastic processes}, volume 131 of
  \emph{Mathematical Surveys and Monographs}.
\newblock American Mathematical Society, Providence, RI, 2006.
\newblock ISBN 978-0-8218-4145-7; 0-8218-4145-9.

\bibitem[Feng et~al.(2012)Feng, Fouque, Kumar, et~al.]{FFK12}
Jin Feng, Jean-Pierre Fouque, Rohini Kumar, et~al.
\newblock Small-time asymptotics for fast mean-reverting stochastic volatility
  models.
\newblock \emph{The Annals of Applied Probability}, 22\penalty0 (4):\penalty0
  1541--1575, 2012.

\bibitem[Fleming and Soner(2006)]{FS06}
Fleming, Wendell H., and Halil Mete Soner. 
\newblock \emph{Controlled Markov processes and viscosity solutions}. Vol. 25. 
\newblock Springer Science \& Business Media, 2006.

\bibitem[Kang and Kurtz(2013)]{KK13}
Hye-Won Kang and Thomas~G. Kurtz.
\newblock Separation of time-scales and model reduction for stochastic reaction
  networks.
\newblock \emph{Ann. Appl. Probab.}, 23\penalty0 (2):\penalty0 529--583, 2013.
\newblock \doi{10.1214/12-AAP841}.
\newblock URL \url{http://projecteuclid.org/euclid.aoap/1360682022}.

\bibitem[Kang et~al.(2014)Kang, Kurtz, and Popovic]{KKP14}
Hye-Won Kang, Thomas~G. Kurtz, and Lea Popovic.
\newblock Central limit theorems and diffusion approximations for multiscale
  {M}arkov chain models.
\newblock \emph{Ann. Appl. Probab.}, 24\penalty0 (2):\penalty0 721--759, 2014.
\newblock ISSN 1050-5164.
\newblock \doi{10.1214/13-AAP934}.
\newblock URL \url{http://dx.doi.org/10.1214/13-AAP934}.

\bibitem[Kontoyannis and Meyn(2005)]{KM05}
Ioannis Kontoyiannis and Sean P. Meyn.
\newblock Large deviations asymptotics and the spectral theory of multiplicatively regular Markov processes.
\newblock \emph{Electron. J. Probab.} 10.3: 61-123, 2005.

\bibitem[Kumar and Popovic(2017)]{KP17}
Rohini Kumar and Lea Popovic.
\newblock Large deviations for multi-scale jump-diffusion processes.
\newblock \emph{Stochastic Processes and their Applications}, 127\penalty0
  (4):\penalty0 1297--1320, 2017.

\bibitem[Lenhart and Yamada(1991)]{Lenhart} 
Suzanne M. Lenhart and Naoki Yamada.
\newblock PerronÕs method for viscosity solutions associated with piecewise-deterministic processes.
\newblock \emph{Funkcialaj Ekvacioj} 34.173-186: 143, 1991.

\bibitem[Li and Lin(2017)]{LL17}
Tiejun Li and Feng Lin.
\newblock Large deviations for two scale chemical kinetic processes.
\newblock \emph{arXiv preprint arXiv:1504.03781}, 2017.

\bibitem[Wentzell(1977)]{Wen77a}
Alexander~D. Wentzell.
\newblock Rough limit theorems on large deviations for {M}arkov stochastic
  processes. {I}.
\newblock \emph{Theory Probab. Appl.}, 21\penalty0 (2):\penalty0 227--242,
  1977.
\newblock \doi{10.1137/1121030}.
\newblock URL \url{http://dx.doi.org/10.1137/1121030}.

\end{thebibliography}

\end{document}